\documentclass{amsart}

\usepackage[frame,cmtip,curve,arrow,matrix,line,graph]{xy}
\usepackage{amssymb}
\usepackage{amscd}
\usepackage{amsmath}

\numberwithin{equation}{section}

\newtheorem{theorem}{Theorem}
\newtheorem{proposition}[theorem]{Proposition}
\newtheorem{corollary}[theorem]{Corollary}
\newtheorem{lemma}[theorem]{Lemma}

\newtheorem{deflem}[theorem]{Definition/Lemma}

\newtheorem{definition}[theorem]{Definition}
\newtheorem{exam}[theorem]{Example}

\newtheorem{rem}[theorem]{Remark}

\newenvironment{remark}{\begin{rem}\rm}{\end{rem}}

\def\lab#1{\label{#1} }
\numberwithin{theorem}{section}

\newcommand{\hatm}{\hat{M}}
\newcommand{\pp}{\mathbb{P}}

\newcommand{\cc}{\mathbb{C}}
\newcommand{\rr}{\mathbb{R}}

\newcommand{\lik}{\mathfrak{k}}

\newcommand{\lih}{\mathfrak{h}}

\newcommand{\codim}{\mathrm{codim \,}}

\newcommand{\git}{/\!\!/}

\newcommand{\Cal}{\mathcal}

\def\cN{\mathcal{N} }

\begin{document}

\title{Intersection cohomology of quotients of nonsingular varieties}
\author{Young-Hoon Kiem}
\address{kiem@math.snu.ac.kr; Department of Mathematics,
Seoul National University, Seoul, 151-747, Korea }

\date{March 2003}
\maketitle

\section{Introduction}
Let $M\subset \pp^n$ be a nonsingular projective variety acted on
by a connected complex reductive group $G=K^{\cc}$ via a
homomorphism $G\to GL(n+1)$ which is the complexification of a
homomorphism $K\to U(n+1)$. Geometric invariant theory (GIT) gives
us a recipe to form a quotient $\phi:M^{ss} \to M\git G$ of the
set of semistable points and many interesting spaces in algebraic
geometry are constructed in this manner \cite{mfk}. But often this
quotient is singular and hence intersection cohomology with middle
perversity is an important topological invariant. The purpose of
this paper is to present a way to compute the middle perversity
intersection cohomology of the singular quotients.

The choice of an embedding $M\subset \pp^n$ provides us with
a moment map for the $K$-action and the Morse stratification
of $M$ with respect to its norm square is $K$-equivariantly
perfect, with the unique open dense stratum $M^{ss}$
 \cite{k1}. Hence one can compute the Betti numbers
for the equivariant cohomology $H^*_K(M^{ss})$ as well as the cup
product, at least in principle. See also \cite{k6}.

When $G$ acts locally freely on $M^{ss}$, we get an orbifold
$M\git G=M^{ss}/G$ and $$H^*_K(M^{ss})\cong H^*(M\git G) \cong
IH^*(M\git G).$$ Hence this equivariant Morse theory enables us to
compute the cohomology ring of the orbifold quotient.

However, if the $G$ action on $M^{ss}$ is not locally free, then
$M\git G$ has more serious singularities than finite quotient
singularities. In general, the induced homomorphism $H^*(M\git
G)\to H^*_K(M^{ss})$, which comes from the quotient map
$${\phi}_K:EK\times_K M^{ss} \to M\git G,$$ is neither injective nor
surjective. Furthermore, the natural map from the ordinary
cohomology $H^*(M\git G)$ to the intersection cohomology
$IH^*(M\git G)$ due to Goresky and MacPherson \cite{GM1,GM2} is
neither injective nor surjective. Hence, knowledge of the
equivariant cohomology $H^*_K(M^{ss})$ does not directly enable us
to compute the topological invariants for $M\git G$.

In \cite{k2}, Kirwan invented a method to partially desingularize
$M\git G$ by blowing up $M^{ss}$ systematically.
When the set of stable points $M^s$ is nonempty, she used this process
to define a map, which we call the Kirwan map
$$\kappa_M^{ss}:H^*_K(M^{ss})\to IH^*(M\git G)$$
and then to provide an algorithm for the computation of the
Betti numbers of $IH^*(M\git G)$. See \cite{k5} for the Betti
number computation of the moduli space of rank 2 holomorphic vector
bundles over a Riemann surface. However, this method does not
give us any information on the intersection pairing of $IH^*(M\git G)$
which is also an essential topological invariant.

In this paper, we construct a natural injection
\begin{equation}\lab{eq1.1}
\phi^*_K:IH^*(M\git G)\to H^*_K(M^{ss})\end{equation} under an
assumption, named the \emph{almost balanced} condition (Definition
\ref{def5.1k}). This homomorphism is the right inverse of the
Kirwan map above (Proposition \ref{prop6.2b}).

It is well-known that intersection cohomology is contravariant
only for some limited classes of maps. One of the most general
known conditions for a subanalytic map
$$f:X\to Y$$
between subanalytic pseudo-manifolds to induce a contravariant
homomorphism
$$f^*:IH^*(Y)\to IH^*(X)$$
by pulling back cycles is the \emph{placid} condition introduced
in \cite{gm6}: The map $f$ is placid if $Y$ has a stratification
such that for each stratum $S$ of $Y$, we have
$$\codim\, S\le \codim\, f^{-1}(S).$$
With this condition, the pull-back of an allowable cycle is again
allowable and hence we get the desired homomorphism $f^*$.

Furthermore, when there is a compact Lie group $K$ acting on $X$
and $f$ is invariant, we have a homomorphism
$$f^*_K:IH^*(Y)\to IH^*_K(X)$$
under the placid condition because $EK\times_K X$ can be
approximated by a sequence of finite dimensional pseudo-manifolds
(see \S\ref{s2.5}).

More generally, if we have
\begin{equation}\lab{pqplain}
q(\codim\,S)\le p(\codim\,f^{-1}(S))
\end{equation}
for some perversities $p$ and $q$, then the pull-back of cycles
gives us a homomorphism
$$f^*:IH^*_q(Y)\to IH^*_p(X).$$
We call (\ref{pqplain}) the $(p,q)$-placid condition.

In our case, $M^{ss}$ is smooth and hence $H^*_K(M^{ss})$ is
isomorphic to the equivariant intersection cohomology
$$IH^*_{p,K}(M^{ss})$$
for any perversity $p$. Hence if $\phi$ is $(p,m)$-placid, we get
a homomorphism
$$\phi^*_K:IH^*(M\git G)\to H^*_K(M^{ss}).$$
Obviously the most general condition is obtained when $p$ is the
top perversity $t$. Hence, when $\phi$ is $(t,m)$-placid,  we
obtain the natural homomorphism (\ref{eq1.1}).

Because $\phi^*_K$ is defined by pulling back cycles, we can
deduce that the intersection pairing is preserved in the following
sense: for $\alpha, \beta$ of complementary degrees in $IH^*(M\git
G)$, i.e. $\deg \alpha+\deg \beta=\dim M\git G$, we have
\begin{equation}\lab{eq1.2}
\phi^*_K(\alpha)\cup \phi^*_K(\beta)=\langle \alpha,\beta\rangle
\phi^*_K(\tau)\end{equation}
 where $\tau$ is the top degree class represented by a point.
Hence, we can compute the intersection numbers in terms of the cup
product structure of the equivariant cohomology.

By the local model theorem (Lemma \ref{local} and Proposition
\ref{localX}) from \cite{SL}, we can describe a topological
stratification of $M^{ss}$ and $M\git G$ explicitly and compute
the codimensions of strata in terms of the weight distributions of
the actions of the stabilizer subgroups on the symplectic slices.
By the results of \cite{k1}, we will see in \S \ref{almost} that
the $(t,m)$-placid condition is in fact a condition on the
balancedness of the weight distributions. This is our almost
balanced condition.

Suppose now the almost balanced condition is satisfied. Then by
(\ref{eq1.1}) and (\ref{eq1.2}), the middle perversity
intersection cohomology of $M\git G$ is completely determined as a
graded vector space with non-degenerate intersection pairing if we
can identify the image of $\phi^*_K$ in $H^*_K(M^{ss})$.

In the light of Deligne's construction of intersection cohomology
sheaf, it seems reasonable to expect that the answer should be
obtained by ``truncating locally" along each stratum. But
complications arise when various strata intersect. To control the
complications, we require that the submanifolds of $M$ fixed by
subgroups of $G$ also satisfy the almost balanced condition. We
call it the \emph{weakly balanced} condition (Definition
\ref{firstWB}). When this condition is satisfied, we can identify
the image of $\phi^*_K$ as a subspace $V^*_M$ of $H^*_K(M^{ss})$
obtained by truncation. See Definition \ref{def7v} and
(\ref{defVker}). By computing the dimension of $V^*_M$, we get the
intersection Betti numbers without going through the
desingularization process as in \cite{k3,k5}. One can furthermore
compute the Hodge numbers since both $\phi^*_K$ and
$\kappa_M^{ss}$ preserve the Hodge structure.

The weakly balanced condition is satisfied by many interesting
spaces including the moduli spaces of holomorphic vector bundles
over a Riemann surface for any rank and degree. (See Proposition
\ref{modspex}.) Also we demonstrate the computations of $V^*_M$
for some standard examples in \S\ref{section9}.

The layout of this paper is as follows. In section 2, we recall
the sheaf theoretic definition of intersection cohomology while in
section 3 we show that intersection cohomology is functorial with
respect to $(p,q)$-placid maps. The stratification of a symplectic
reduction is discussed in section 4 and the homomorphism
$\phi^*_K$ is constructed in section 5 after defining the almost
balanced condition. In section 6 we recall the Kirwan map and
prove that $\phi^*_K$ is the right inverse of $\kappa_M^{ss}$. In
section 7 we define $V_M^*$ and state the theorem, which says
\begin{equation}\label{eq1.3}\phi^*_K(IH^*(M\git
G))=V_M^*.\end{equation} This is proved in section 8 and several
examples are computed in section 9.

The proof of (\ref{eq1.3}) is unfortunately quite lengthy and so I
briefly sketch the outline. Our proof is by induction on the
maximum among the dimensions of stabilizers. We first show that
$$\phi^*_K(IH^*(M\git G))\subset V_M^*.$$
This is best seen from the sheaf theoretic perspective and we need
the full strength of the weakly balanced condition for the
computation of stalks.

Next let $\hatm$ be the first blow-up in the partial
desingularization process \cite{k2}. By our induction hypothesis,
we have $$\hat{\phi}^*_K(IH^*(\hatm\git G))=V^*_{\hatm}$$ where
$\hat{\phi}:\hatm^{ss}\to \hatm\git G$ is the GIT quotient map.
Then we show there is an embedding
$$V^*_M\hookrightarrow V^*_{\hatm}\cong IH^*(\hatm\git G).$$
Since the difference between $IH^*(\hatm\git G)$ and $IH^*(M\git
G)$ comes from a sheaf complex supported at the exceptional
divisor, we will see that the proof reduces to the computation for
the normal bundle $\cN$ to the blow-up center in $M^{ss}$. By
spectral sequence, this further reduces to the computation for the
normal space at a point. Upon checking this, we complete the
proof.

The main application of the results of this paper is to generalize
the theorem of Jeffrey and Kirwan \cite{JK98} about the
intersection numbers on the \emph{smooth} moduli spaces of bundles
over a Riemann surface. In a joint work with Jeffrey, Kirwan and
Woolf, we will show that the intersection numbers of the
intersection cohomology of the \emph{singular} moduli spaces of
bundles over a Riemann surface are given by a residue formula
similar to that of \cite{JK98}.

The first version of this paper was written several years ago.
There $V^*_M$ was defined and it was shown that the Kirwan map
$\kappa_M^{ss}$ restricted to $V^*_M$ is an isomorphism onto
$IH^*(M\git G)$ under the weakly balanced condition (Definition
7.2). About a year later, in \cite{kw}, we generalized the Kirwan
map to symplectic reductions and extended the results of this
paper to the purely symplectic setting by interpreting the first
condition of the weakly balanced action as the cosupport axiom for
intersection homology sheaf. After finishing \cite{kw}, it
occurred to us that the almost-balanced condition could be most
naturally explained in terms of $(p,q)$-placid maps and hence the
current paper was revised using this new observation. In
\cite{kiem2}, we showed that the extended moduli spaces
 defined by L. Jeffrey \cite{j} satisfy our assumptions and
hence we can compute the intersection cohomology of representation spaces
of surface groups in terms of the equivariant
cohomology.

Every cohomology group in this paper has \emph{complex
coefficients}.

{\bf Acknowledgements.} I would like to express my deep gratitude
to Professor Ronnie Lee for encouragement and advice. I am
grateful to Professors Frances Kirwan, Lisa Jeffrey and Jon Woolf
for useful discussions and crucial comments. The referee's
suggestions improved the exposition considerably and I wish to
thank the referee.


\section{Intersection cohomology}
Intersection cohomology was introduced by Goresky and MacPherson
in \cite{GM1, GM2} as an invariant of singular spaces which
retains useful properties like Poincar\'e duality and Lefschetz
theorems. In this section we recall the definition and some
properties that we will use.

\subsection{Topological stratification}\lab{s2.1}
 An even dimensional topological space $X$ equipped with a
filtration by even dimensional closed subsets
\begin{equation}\lab{eq2.1}(\mathfrak{X}) \ \ \ \
X=X_n\supset X_{n-2}\supset  \cdots\supset X_0\supset
\emptyset\end{equation} is called a \emph{stratified
pseudo-manifold} if
\begin{itemize}\item $X-X_{n-2}$ is dense
\item  for each $i$, $S_i=X_i- X_{i-2}$ is a topological manifold of dimension $i$  or empty
\item for each $x\in S_i$ there are a compact stratified
pseudo-manifold $$L=L_{n-i-1}\supset \cdots\supset L_0\supset
\emptyset$$ and a stratum-preserving homeomorphism of a
neighborhood of $x$ onto $\rr^i\times cL$ where $cL=L\times
[0,\infty)/L\times 0$.
\end{itemize}

In this case, $\mathfrak X$ is called a \emph{topological
stratification}. A pseudo-manifold $X$ of dimension $n$ is normal
if $H_n(X,X-x)=1$ for any $x\in X$.

\subsection{Deligne's construction of intersection cohomology sheaf}\lab{s2.2}
Let
$$p:\{2,3,\cdots, n\}\to \{0,1,2,\cdots,n-2\}$$ be an increasing function such
that $q(i)=i-2-p(i)$ is also an increasing nonnegative function.
We call such $p$ a perversity. For instance,
$$m(i)=[\frac{i-2}2]$$ is the middle perversity and $t(i)=i-2$ is the top perversity.
Given a normal stratified pseudo-manifold $X$ with stratification
$(\mathfrak{X})$, put $U_{2k}=X- X_{n-2k}$ and let
$j_{2k}:U_{2k}\hookrightarrow U_{2k+2}$ denote the inclusion. Then
the perversity $p$ intersection cohomology $IH^*_p(X)$ of $X$ is
the hypercohomology of the sheaf complex defined inductively by
$$\mathcal{IC}^{\cdot}_{p,X}|_{U_{2}}=\cc_{U_2}$$
$$\mathcal{IC}^{\cdot}_{p,X}|_{U_{2k+2}}\cong \tau_{\le p(2k)}R{j_{2k}}_*
(\mathcal{IC}^{\cdot}_{p,X}|_{U_{2k}}).$$ This defines an object
in the derived category $\bold{D}^+_c(X)$ of bounded below
cohomologically constructible sheaves, which is independent of the
choice of stratification, and the intersection cohomology is a
homeomorphism invariant \cite{GM2}. If $p$ is the zero perversity,
i.e. $p(i)=0$, then $\mathcal{IC}^{\cdot}_{p,X}\cong \cc$ and
hence $IH^*_0(X)$ is isomorphic to the (singular) cohomology
$H^*(X)$. When $p$ is the top perversity $t$,
$\mathcal{IC}^{\cdot}_{t,X}$ is isomorphic to the dualizing
complex $D_X^{\cdot}$ and thus $IH^i_t(X)$ is isomorphic to the
Borel-Moore homology $H^{BM}_{n-i}(X)$.

More generally, suppose we have a filtration
$$(\mathfrak{X}')\ \ \ \ X=X_n'\supset X'_{n-2}\supset \cdots\supset X'_0\supset \emptyset $$
by even dimensional closed subsets (not necessarily a topological
stratification) such that \begin{enumerate}
\item $X-X'_{n-2}$ is dense,
\item $S'_i=X'_i-X'_{i-2}$ is a $\cc$-homology manifold of
dimension $i$ or empty.
\end{enumerate}

Let $U'_{2k}=X- X'_{n-2k}$ and $j'_{2k}:U'_{2k}\hookrightarrow
U'_{2k+2}$. Define a sheaf complex $\mathcal{P}^{\cdot}_{p,X}$
inductively by
$$\mathcal{P}^{\cdot}_{p,X}|_{U'_{2}}=\cc_{U'_2}$$
$$\mathcal{P}^{\cdot}_{p,X}|_{U'_{2k+2}}\cong \tau_{\le p(2k)}R{j'_{2k}}_*
(\mathcal{P}^{\cdot}_{p,X}|_{U'_{2k}}).$$

\begin{deflem}
If $\mathcal{P}^{\cdot}_{p,X}$ is topologically constructible (see
\cite{GM2} p83) then $\mathcal{P}^{\cdot}_{p,X}\cong
\mathcal{IC}^{\cdot}_{p,X}$. We call $\mathfrak{X}'$ a \emph{nice}
filtration if $\mathcal{P}^{\cdot}_{p,X}$ is topologically
constructible.
\end{deflem}
\begin{proof}
It suffices to show that $\mathcal{P}^{\cdot}_{p,X}$ satisfies
[AX2] in \cite{GM2} p107. The proof is exactly same as the proof
of Lemma 2 in \cite{GM2} p110.
\end{proof}

For a space $X$ constructed as the GIT quotient of a smooth
variety, the decomposition by orbit types gives us a topological
stratification while the infinitesimal orbit types will give us
only a nice filtration. See \S\ref{section4}.

\subsection{Geometric chains}\lab{s2.3}
Explicitly intersection cohomology can by described in terms of
geometric chains. Suppose $X$ is a subanalytic pseudo-manifold
with subanalytic stratification $\mathfrak{X}=\{X_i\}$. For an
open subset $U$ let $C^*(U)$ be the chain complex defined by
$$C^i(U)=\{\text{locally finite chains in }U\text{ of dimension }n-i \}.$$
Note that our indexing scheme is ``cohomology superscript"
(\cite{GM2} p98).

For an open subset $U$ and a perversity $p$, define a subcomplex
of $C^*(U)$ by
$$\begin{array}{ll}
IC^i_p(U)=\{\xi\in C^i(U)\,|\,&\dim(|\xi|\cap X_{n-c}\cap U)\le
n-i-c+p(c),\\ &\dim(|\partial \xi|\cap X_{n-c}\cap U)\le
n-i-1-c+p(c)\}.\end{array}$$ This gives rise to a \emph{soft}
sheaf complex which is isomorphic to $\mathcal{IC}^{\cdot}_{p,X}$.
(See \S2.1, \S3.6 \cite{GM2}.)
 Hence the cohomology of $IC^*_p(X)$ is
the perversity $p$ intersection cohomology $IH^*_p(X)$.

\subsection{Intersection pairing}\lab{s2.4}
If perversities $p,q,r$ satisfy $p+q\le r$, then there is a unique
morphism $$\mathcal{IC}^{\cdot}_{p,X}\otimes^L
\mathcal{IC}^{\cdot}_{q,X}\to \mathcal{IC}^{\cdot}_{r,X}$$ in
$\bold{D}^+_c(X)$ extending $\cc\otimes \cc\to \cc$ over $U_2$
(\cite{GM2} p112). This morphism induces a homomorphism
$$IH^i_p(X)\otimes IH^j_q(X)\to IH^{i+j}_r(X).$$ For instance if
$p=q=r=0$ then this is just the cup product of cohomology classes.

For the middle perversity $m$, we have $m+m\le t$. Thus if $X$ is
normal compact connected oriented pseudo-manifold, we have the
intersection pairing
\begin{equation} IH^i_m(X)\otimes
IH^j_m(X)\to IH^n_t(X)\cong \cc
\end{equation}
for $i+j=n$. In terms of geometric chains, this intersection
pairing is defined as follows: For $\alpha, \beta\in IH^*_m(X)$
such that $\deg(\alpha)+\deg(\beta)=n$, we can find representative
cycles $\xi$ and $\sigma$ such that they intersect only at
finitely many points in the smooth part $U_2$, transversely. The
intersection pairing $\langle\alpha,\beta\rangle$ of $\alpha$ and
$\beta$ is the number of intersection points counted with signs as
usual. (See \S2.3 \cite{GM1}.)

\subsection{Equivariant intersection cohomology}\lab{s2.5}
Suppose a compact connected Lie group $K$ acts on a
pseudo-manifold preserving a topological stratification
$\mathfrak{X}=\{X_i\}$ (\S\ref{s2.2}). Let $EK$ be a contractible
space on which $K$ acts freely and $BK=EK/K$. Then the closed
subsets
$$EK\times_KX_i$$ form a filtration of $X_K:=EK\times_KX$. For a perversity $p$,
apply Deligne's construction (\S\ref{s2.2}) to this filtration to
obtain a sheaf complex $\mathcal{IC}^{\cdot}_{p,X_K}$. The
equivariant intersection cohomology $IH^*_{p,K}(X)$ is defined as
the hypercohomology of this sheaf complex. (See \S5.2, \S13.4 in
\cite{BL}.) From the fibration $X_K\to BK$ with fiber $X$, we get
a spectral sequence
$$H^i(BK)\otimes IH^j(X)\Rightarrow IH^{i+j}_{p,K}(X).$$

More concretely, choose a smooth classifying sequence (\S12
\cite{BL})
$$EK_0\subset EK_1\subset EK_2\subset \cdots$$ where $EK_k$ is a
$k$-acyclic free $K$-manifold. Let $BK_k=EK_k/K$ so that
$BK=\lim_{k\to\infty}BK_k$. The filtration $$EK_k\times_KX_i$$ of
$EK_k\times_KX$ is a topological stratification and thus Deligne's
construction gives rise to
$\mathcal{IC}^{\cdot}_{p,EK_k\times_KX}$. From the spectral
sequences for the fibrations
$$\xymatrix{
X\ar[r]&EK_k\times_KX\ar[d]\\
&BK_k }\qquad \xymatrix{
X\ar[r]&EK\times_KX\ar[d]\\
&BK }$$ we see immediately that
\begin{equation}\lab{eq2l}
IH^*_{p,K}(X)=\lim_{k\to\infty}IH^*_p(EK_k\times_KX)\end{equation}
because $H^{<k}(BK_k)\cong H^{<k}(BK)$ by Whitehead's theorem.


\section{Placid maps}
Unlike ordinary cohomology, intersection cohomology is
contravariant only for a limited class of maps. In this section we
generalize slightly the concept of placid maps due to Goresky and
MacPherson and show the functoriality. This will give us the right
inverse of the Kirwan map.

\subsection{Placid maps}\lab{s3.1}
Let $f:X\to Y$ be a subanalytic map between subanalytic
pseudo-manifolds. Suppose $Y$ is compact.

\begin{definition}
$f$ is called \emph{$(p,q)$-placid} if there is a subanalytic
stratification of $Y$ such that
\begin{equation}q(\codim S)\le p(\codim f^{-1}(S))\lab{eq3.1}\end{equation}
 for each stratum $S$ of $Y$.
\end{definition}

When $p=q$, we recover the original placid maps \cite{gm6}.

With this definition, Proposition 4.1 in \cite{gm6} is modified as
follows.

\begin{proposition}\lab{funct}
If $f$ is $(p,q)$-placid, then the pull-back of generic chains
induces a homomorphism on intersection cohomology
\begin{equation}f^*:IH^*_q(Y)\to IH^*_p(X).\lab{eq3.2}\end{equation}
\end{proposition}
\begin{proof}
See p373 in \cite{gm6} for details. Since $f$ is subanalytic,
there is a stratification of $X$ for which $f$ is a stratified map
(\cite{GM2} \S1.2). By McCrory's transversality, any cohomology
class $\alpha$ in $IH^i_q(Y)$ can be represented by a chain $\xi$
which is dimensionally transverse to any stratum in $X$ and the
cycle $f^{-1}(\xi)$ lies in $IC^i_p(X)$ because of (\ref{eq3.1}).
The class $f^*\alpha$ is represented by the cycle $f^{-1}(\xi)$.
\end{proof}

\subsection{Pull-back morphism}\lab{s3.2}
The homomorphism (\ref{eq3.2}) comes from a morphism
\begin{equation}\lab{eq3.3}
f^*:\mathcal{IC}^{\cdot}_{q,Y}\to
Rf_*\mathcal{IC}^{\cdot}_{p,X}\end{equation} in the derived
category $\bold{D}^+(Y)$.

Suppose $\mathfrak{X}=\{X_i\}_{i=0}^n$ (resp.
$\mathfrak{Y}=\{Y_i\}_{i=0}^l$) is a \emph{nice} filtration of $X$
(resp. $Y$). Let $f:X\to Y$ be a continuous map such that for
every connected component $S$ of any stratum $Y_i-Y_{i-2}$,
$f^{-1}(S)$ is a union of some connected components of strata in
$X$ and (\ref{eq3.1}) holds.

Over $U_2=Y-Y_{l-2}$, let $f^*_{U_2}$ be the composition of the
adjunction morphism $\cc_{U_2}\to
Rf_*f^*\cc_{U_2}=Rf_*\cc_{f^{-1}(U_2)}$ with the morphism
$Rf_*\cc_{f^{-1}(U_2)}\to
Rf_*\mathcal{IC}^{\cdot}_{p,X}|_{f^{-1}(U_2)}$ from \S5.1
\cite{GM2}.

\begin{proposition} \lab{p3.2} There is a unique morphism (\ref{eq3.3})
in $\bold{D}^+_c(Y)$ extending $f^*_{U_2}$.
\end{proposition}
We need the following simple lemma to prove the proposition.

\begin{lemma}\lab{l3.1}
For $\mathcal{A}^{\cdot}\in \bold{D}^+_c(X)$, $\tau_{\le
q}Rf_*\mathcal{A}^{\cdot}\cong \tau_{\le q}Rf_*\tau_{\le
p}\mathcal{A}^{\cdot}$ if $q\le p$.
\end{lemma}
\begin{proof}
Choose an injective resolution $\mathcal{I}^{\cdot}$ of
$\mathcal{A}^{\cdot}$. (See for instance \cite{GeM} p181.) It is
elementary to find an injective resolution $\mathcal{J}^{\cdot}$
of $\tau_{\le p}\mathcal{A}^{\cdot}$ such that
$\mathcal{I}^{k}=\mathcal{J}^{k}$ for $k\le p$. Then
$Rf_*\mathcal{A}^{\cdot}=f_*\mathcal{I}^{\cdot}$ is equal to
$Rf_*\tau_{\le p}\mathcal{A}^{\cdot}=f_*\mathcal{J}^{\cdot}$ up to
degree $p$. So we proved the lemma.\end{proof}

\begin{proof}[Proof of Proposition \ref{p3.2}]
Let $U_{2k}=Y-Y_{l-2k}$ and $j_{2k}:U_{2k}\hookrightarrow
U_{2k+2}$. Recall from \S2.2 that Deligne's construction applied
to $\mathfrak Y$ with perversity $q$ gives us the intersection
cohomology sheaf $\mathcal{IC}_{q,Y}^{\cdot}$. Suppose we have
constructed
$$f^*_{U_{2k}}:\mathcal{IC}^{\cdot}_{q,Y}|_{U_{2k}}\to
Rf_*\mathcal{IC}^{\cdot}_{p,X}|_{f^{-1}(U_{2k})}.$$ This gives
rise to
$$\begin{array}{ll}
\mathcal{IC}^{\cdot}_{q,Y}|_{U_{2k+2}}&\cong \tau_{\le
q(2k)}{Rj_{2k}} _*\mathcal{IC}^{\cdot}_{q,Y}|_{U_{2k}}\to
\tau_{\le
q(2k)}{Rj_{2k}} _*Rf_*\mathcal{IC}^{\cdot}_{p,X}|_{f^{-1}(U_{2k})}\\
&=\tau_{\le q(2k)}Rf_*{Ri_{2k}}
_*\mathcal{IC}^{\cdot}_{p,X}|_{f^{-1}(U_{2k})}\\
&\cong \tau_{\le q(2k)}Rf_*\tau_{\le q(2k)}{Ri_{2k}}
_*\mathcal{IC}^{\cdot}_{p,X}|_{f^{-1}(U_{2k})}\end{array}$$ by
Lemma \ref{l3.1} where $i_{2k}:f^{-1}(U_{2k})\hookrightarrow
f^{-1}(U_{2k+2})$ since $j_{2k}\circ f=f\circ i_{2k}$. We claim
$$\tau_{\le q(2k)}{Ri_{2k}}
_*\mathcal{IC}^{\cdot}_{p,X}|_{f^{-1}(U_{2k})}\cong \tau_{\le
q(2k)} \mathcal{IC}^{\cdot}_{p,X}|_{f^{-1}(U_{2k+2})}.$$

For simplicity, suppose $f^{-1}(U_{2k+2})-f^{-1}(U_{2k})$ consists
of only one connected stratum $\tilde S$. Then
$$\mathcal{IC}_{p,X}|_{f^{-1}(U_{2k+2})}\cong \tau^{\tilde S}_{\le p(\codim
\tilde{S})}R{i_{2k}}_*\mathcal{IC}_{p,X}|_{f^{-1}(U_{2k})}
$$ where $\tau^{\tilde S}_{\le p}$ is the ``truncation over a closed subset
functor" (see \S1.14 \cite{GM2}). In particular,
$$\begin{array}{ll}
\tau_{\le q(2k)}\mathcal{IC}_{p,X}|_{f^{-1}(U_{2k+2})}&\cong
\tau_{\le q(2k)}\tau^{\tilde S}_{\le p(\codim
\tilde{S})}R{i_{2k}}_*\mathcal{IC}_{p,X}|_{f^{-1}(U_{2k})}\\
&\cong \tau_{\le
q(2k)}R{i_{2k}}_*\mathcal{IC}_{p,X}|_{f^{-1}(U_{2k})}\end{array}
$$ because $q(2k)\le p(\codim {\tilde S})$.  When there are more than one strata in
$f^{-1}(U_{2k+2})-f^{-1}(U_{2k})$ we simply repeat the argument
for each stratum in the order of increasing codimension. Hence we
get a morphism
$$\begin{array}{ll}
\mathcal{IC}^{\cdot}_{q,Y}|_{U_{2k+2}}&\to \tau_{\le q(2k)}Rf_*
\tau_{\le q(2k)} \mathcal{IC}^{\cdot}_{p,X}|_{f^{-1}(U_{2k+2})}\\
& \cong \tau_{\le q(2k)}Rf_*
\mathcal{IC}^{\cdot}_{p,X}|_{f^{-1}(U_{2k+2})}\to Rf_*
\mathcal{IC}^{\cdot}_{p,X}|_{f^{-1}(U_{2k+2})}.
\end{array}$$ Uniqueness is an elementary exercise.
\end{proof}

By taking hypercohomology, the morphism (\ref{eq3.3}) induces a
homomorphism
$$IH_q^*(Y)\to IH_p^*(X).$$
Suppose $f$ is subanalytic and $\mathfrak Y$ is a subanalytic
stratification. Note that an intersection cycle $\xi$ is
completely determined by its intersection with the open dense
stratum (\cite{bo} p10) and over $U_2$ (\ref{eq3.3}) is just the
adjunction. If we use the sheaf complexes by geometric chains
(\S\ref{s2.3}) for $\mathcal{IC}^{\cdot}_{p,X}$ and
$\mathcal{IC}^{\cdot}_{q,Y}$, the induced homomorphism sends the
class $[\xi]\in IH_q^*(Y)$, which is represented by a chain $\xi$
dimensionally transverse to any stratum of $X$, to
$[f^{-1}(\xi)]$. In other words, the induced homomorphism is
exactly the homomorphism (\ref{eq3.2}).

\subsection{Equivariant case}
Suppose a compact connected Lie group $K$ acts on $X$ preserving a
nice filtration $\mathfrak{X}=\{X_i\}_{i=0}^n$. Deligne's
construction with the filtration $\{EK\times_KX_i\}$ and
perversity $p$ gives us $\mathcal{IC}^{\cdot}_{p,X_K}$ in
$\bold{D}^+(X_K)$ whose hypercohomology is $IH^*_{p,K}(X)$.

Let $Y$ be a compact pseudo-manifold with a nice filtration
$\mathfrak{Y}=\{Y_i\}_{i=0}^m$. Let $f:X\to Y$ be an invariant
continuous map such that for any connected component $S$ of a
stratum, $f^{-1}(S)$ is a union of some connected components of
strata in $X$ and (\ref{eq3.1}) holds.

Then the proof of Proposition \ref{p3.2} gives us a morphism
\begin{equation}\lab{eq3.4a}\mathcal{IC}_{q,Y}^{\cdot} \to
R{f_K}_*\mathcal{IC}_{p,X_K}^{\cdot}\end{equation} where
$f_K:X_K=EK\times _K X \to Y$ is the obvious map induced from $f$.
So we have a homomorphism
\begin{equation}\lab{eq3.5a}f^*_K:IH^*_q(Y)\to IH^*_{p,K}(X).
\end{equation}

Similarly, we have a morphism
$$\mathcal{IC}_{q,Y}^{\cdot}
\to R{f_k}_*\mathcal{IC}_{p,EK_k\times_KX}^{\cdot}$$ where
$f_k:EK_k\times _K X \to Y$ and a homomorphism
$$f^*_k:IH^*_q(Y)\to IH^*_{p}(EK_k\times_KX).
$$
From the commutative diagram $$\xymatrix{
EK_k\times_KX\ar[r]^{\imath_k}\ar[dr]&EK\times_KX\ar[d]\\
&Y }$$ we see that the composition
$$\xymatrix{IH^*_q(Y)\ar[r]^{f^*_K}& IH^*_{p,K}(X)\ar[r]^(.4){\imath_k^*}
&IH^*_p(EK_k\times_KX)}$$ is $f^*_k$ since
$\imath_k^*\mathcal{IC}_{p,X_K}^{\cdot}\cong
\mathcal{IC}_{p,EK_k\times_KX}^{\cdot}$ by construction. Hence
$\{f_k^*\}$ determine $f^*_K$ because of (\ref{eq2l}).

\subsection{$f^*$ preserves the intersection pairing}\lab{s3.3}
Our interest lies in the case where $X$ is a smooth analytic
manifold with an action by a compact connected Lie group $K$ and
$f:X\to Y$ is invariant. In particular we wish to relate the
middle perversity intersection cohomology $IH^*_m(Y)$ with the
equivariant cohomology $H^*_K(X)=H^*(EK\times_KX)$. From now on
when using middle perversity, we will drop the subscript $m$ for
convenience.

Suppose $X$ is smooth. Then since intersection cohomology is
independent of stratification, $IH^*_{p}(X)\cong H^*(X)$ for any
perversity $p$ and hence for any $(p,m)$-placid subanalytic map
$f$ we have a homomorphism $f^*:IH^*_q(Y)\to H^*(X)$. Obviously
$(t,m)$-placid condition is most general for us to get such a
homomorphism.

If furthermore $K$ acts on $X$ preserving a subanalytic
stratification $\mathfrak X$ for which $f$ is $(t,m)$-placid
stratified map with a stratification $\mathfrak Y$ of $Y$, we have
a morphism $$f^*_K:\mathcal{IC}_{Y}^{\cdot} \to
R{f_K}_*\cc_{X_K}$$ and the induced homomorphism
$$f^*_K:IH^*(Y)\to H^*_K(X).$$

\begin{proposition}\lab{p3.5}
 Let $f:X\to Y$ be a $(t,m)$-placid map.
Suppose $X$ is smooth and $Y$ is compact connected normal
oriented. Let $\tau$ be the top degree class represented by a
point in the smooth part. Then the induced map $f^*:IH^*(Y)\to
H^*(X)$ preserves the intersection pairing in the sense that
\begin{equation}\lab{eq3.4} f^*(\alpha)\cup f^*(\beta)=\langle
\alpha,\beta\rangle f^*(\tau)\end{equation} for any
$\alpha,\beta\in IH^*(Y)$ of complementary degrees. In the
equivariant case, the same is true for $f_K^*:IH^*(Y)\to
H^*_K(X)$.
\end{proposition}
\begin{proof}
Recall from \S\ref{s2.4} that $\alpha$ and $\beta$ are represented
by intersection cycles $\xi$ and $\sigma$ that intersect only at
finitely many points in the smooth part. In this case,
$\xi\cap\sigma\in IC^n(Y)$ represents $\langle\alpha,\beta\rangle
\tau\in IH^n(Y)$. By McCrory's transversality result we can
further assume that $\xi$ and $\sigma$ are dimensionally
transverse to each stratum in $X$ so that $f^*\alpha$ and
$f^*\beta$ are represented by $f^{-1}(\xi)$ and $f^{-1}(\sigma)$.
Because $X$ is smooth, the complex $\mathcal{C}^{\cdot}$ of
geometric chains is isomorphic to the constant sheaf $\cc_X$ and
the cup product is just the intersection of chains. Hence the cup
product $f^*(\alpha)\cup f^*(\beta)$ is represented by
$f^{-1}(\xi)\cap f^{-1}(\sigma)=f^{-1}(\xi\cap \sigma)$ which also
represents $\langle \alpha,\beta\rangle f^*(\tau)$. So we proved
(\ref{eq3.4}).

For the equivariant case, observe that the statement is true for
$f_k:EK_k\times_KX\to Y$ for any $k$ since $EK_k\times_KX$ is a
finite dimensional manifold. If we take a sufficiently large $k$,
then $H^{\le l}_K(X)\cong H^{\le l}(EK_k\times_KX)$ and
$f^*_k=f^*_K$ where $l=\dim Y$. So we are done.
\end{proof}

\begin{corollary}\lab{c3.6}
Suppose furthermore $f^*(\tau)\ne 0$ in $H^*(X)$. Then $f^*:IH^*(Y)\to H^*(X)$
is injective and the intersection pairing is given by the cup product structure
of $H^*(X)$. A similar result is true for the equivariant case.
\end{corollary}
\begin{proof}
The result follows from Proposition \ref{p3.5} since the
intersection pairing is non-degenerate for $IH^*(Y)$.\end{proof}


\section{Symplectic reduction}\lab{section4}
Let $(M,\omega)$ be a connected Hamiltonian $K$-space with proper
moment map $\mu:M\to \lik^*$ where $\lik =Lie(K)$. Then the
symplectic reduction $X=\mu^{-1}(0)/K$, which we denote by $M\git
K$, is in general a pseudomanifold, whose strata are symplectic
manifolds. In this section, we describe the orbit type
stratification and the infinitesimal orbit type decomposition of
$X$ from \cite{SL,mes}.

\subsection{Stratification of $X$}
Let $Z=\mu^{-1}(0)$ and $Z_{(H)}=\{x\in Z\, |\,
\mathrm{Stab}(x)\in (H)\}$ for a subgroup $H$ of $K$, where $(H)$
denotes the conjugacy class of $H$. Also, let $Z_H=\{x\in Z\, |\,
\mathrm{Stab}(x)=H\}$. Then
\begin{equation}\lab{eq4.7a}
Z=\bigsqcup_{(H)} Z_{(H)}\end{equation} and
$Z_{(H)}=K\times_{N^H}Z_H$ where $N^H$ is the normalizer of $H$ in
$K$. This decomposition induces a stratification
\begin{equation}\lab{eq4.8a}
X=\bigsqcup_{(H)}X_{(H)}\end{equation}
 where
$X_{(H)}=Z_{(H)}/K=Z_H/N^H$.

For $x\in Z_H$, consider the \emph{symplectic slice}
\begin{equation}\lab{eqsympsl}T_x(Kx)^{\perp_\omega}/T_x(Kx)\end{equation}
 and let $W$ be the symplectic
complement of the $H$-fixed point set in the slice. We recall the
following result from \cite{SL} Lemma 7.1.

\begin{lemma}\lab{local}
There exists a neighborhood of the submanifold $Z_{(H)}$ of $M$ that is symplectically
and $K$-equivariantly diffeomorphic to a neighborhood of the zero section of a vector
bundle $\mathcal{N}\to Z_{(H)}$. The space $\mathcal{N}$ is a symplectic fiber bundle
over the stratum $X_{(H)}$
$$F\to \mathcal{N}\to X_{(H)}$$
with fiber $F$ given by $$F=K\times_H((\lik/\lih)^*\times W).$$
\end{lemma}

In particular, for any $x\in Z_H$, there is a neighborhood of the
orbit $Kx$ that is equivariantly diffeomorphic to $F\times
\rr^{\dim X_{(H)}}$.

Now the reduction of $\mathcal{N}$ is homeomorphic to a
neighborhood  of $X_{(H)}$ in $X$. The principle of reduction in
stages gives us the following.

\begin{proposition}\lab{localX}\textrm{(\cite{SL} 7.4)}
 Given a stratum $X_{(H)}$ of $X$,
there exists a fiber bundle over $X_{(H)}$ with typical fiber
being the cone $W\git H$ such that a neighborhood of the vertex
section of this bundle is symplectically diffeomorphic to a
neighborhood of the stratum inside $X$.
\end{proposition}
Consequently, $\{X_{(H)}\}$ gives us a topological stratification
of $X$. This is called the \emph{orbit type stratification}.

\subsection{Stratification of $M^{ss}$}\lab{subs4.2}
Let $M^{ss}$ denote the open subset of elements in $M$ whose
gradient flow for $f=-|\mu|^2$ has a limit point in $Z$ and put
$r:M^{ss}\to Z$ denote the retraction by the flow.

Let $\phi$ be the composition
\begin{equation}\lab{eq4.9a}
M^{ss}\to Z\to Z/K=X \end{equation} of the retraction $r$ and the
quotient map. Let us call it the symplectic quotient map. The
inverse image $\phi^{-1}(X_{(H)})$ of the stratum $X_{(H)}$ is
diffeomorphic to a sub-fiber bundle of $\mathcal{N}$ in Lemma
\ref{local} with typical fiber $K\times_H((\lik/\lih)^*\times
\phi_W^{-1}(*))$ where $\phi_W:W\to W\git H$ is the symplectic
quotient map for $W$ and $*=\phi_W(0)$ is the vertex of the cone
$W\git H$. If we assign a complex structure, compatible with the
symplectic structure, it is well-known that $\phi_W^{-1}(*)$ is
the affine cone over $\pp W- \pp W^{ss}$ where the superscript
$ss$ denotes the semistable set defined by Mumford \cite{mfk}.
Hence, the affine cones over the unstable strata of $\pp W$ minus
0, together with $\{0\}$, give us a stratification of
$\phi^{-1}(X_{(H)})$ via the diffeomorphism in Lemma \ref{local}.
Observe that $M^{ss}$ is diffeomorphic to a neighborhood of $Z$ by
the gradient flow of $-|\mu|^2$. Since the diffeomorphism in Lemma
\ref{local} is $K$-equivariant and the stratification of
$\phi^{-1}(X_{(H)})$  is completely determined by the group
action, we get a K-invariant stratification of $M^{ss}$ for which
$\phi:M^{ss}\to X$ is stratified.

See \cite{k7} for a description of
the above stratification for GIT quotients and
an application to the Atiyah-Jones conjecture.

\subsection{Infinitesimal orbit type decomposition}
There is another way to decompose $X$ which is useful for partial
desingularization. For a Lie subalgebra $\mathfrak{h}$ of
$\mathfrak{k}$, let $Z_{\mathfrak{h}}=\{x\in
Z\,|\,\text{Lie}\text{Stab}(x)=\mathfrak{h}\}$ and
$$Z_{(\mathfrak{h})}=\{x\in
Z\,|\,\text{Lie}\text{Stab}(x)\in(\mathfrak{h})\}$$ where
$(\mathfrak{h})$ is the conjugacy class of $\mathfrak{h}$. Let
$X_{(\mathfrak{h})}=Z_{(\mathfrak{h})}/K$. Then we have a
decomposition
$$X=\bigsqcup_{(\mathfrak{h})}X_{(\mathfrak{h})}.$$
This is called the infinitesimal orbit type decomposition of $X$.
From \cite{mes} \S3, $X_{(\mathfrak{h})}$ are just orbifolds and
thus homology manifolds. If we use the local normal form in
\cite{mes} \S3.1, it is elementary to show that Deligne's
construction for this decomposition gives us a topologically
constructible sheaf complex with respect to the orbit type
stratification. This fact will not be used in this paper so we
leave the details to the reader.


\section{Almost-balanced action}\lab{almost}

We use the notations of the previous section. We assume that there
is at least one point in $Z=\mu^{-1}(0)$ with finite stabilizer.

\subsection{Placid maps for symplectic quotients}
With the orbit type stratification, $\phi^{-1}(X_{(H)})$ is a
fiber bundle over $X_{(H)}$ with fiber $K\times_H
((\lik/\lih)^*\times \phi_W^{-1}(*))$ and hence
$$\codim_{M^{ss}}\, \phi^{-1}(X_{(H)})=\codim_W\,
\phi_W^{-1}(*).$$ On the other hand,
$$\codim_X\, X_{(H)}=\dim W\git H$$ by Proposition \ref{localX}.
Recall that the map  $\phi:M^{ss}\to X$ is $(t,m)$-placid if
$$t(\codim_{M^{ss}}\, \phi^{-1}(X_{(H)}))\ge m(\codim_X\, X_{(H)})$$
or equivalently
$$\codim_W\, \phi_W^{-1}(*)-2\ge \frac12 \dim W\git H -1.$$
Hence $\phi$ is $(t,m)$-placid if and only if
\begin{equation}\lab{eq5.1c}
\codim_W\, \phi_W^{-1}(*) > \frac12 \dim W\git H\end{equation}
 for each $(H)$.

The unstable strata of $\pp W$ by the norm square of the moment
map can be described using the weights of the maximal torus
action, as follows: For any collection of weights of the maximal
torus action, we consider the convex hull of them and get the
closest point from the origin to the hull. Let $\mathcal{B}$ be
the set of such closest points in the positive Weyl chamber. Then
the unstable strata are in one-to-one correspondence with the set
$\mathcal{B}$. (See \cite{k1}.)

For each $\beta\in \mathcal{B}$, let $n(\beta)$ denote the number
of weights $\alpha$ such that $\langle \alpha,\beta\rangle
<\langle \beta,\beta\rangle$. Then Kirwan proved in \cite{k1} that
the codimension of the stratum corresponding to $\beta\in
\mathcal{B}$ is precisely
$$2n(\beta)-\dim H/\mathrm{Stab}\beta$$
where $\mathrm{Stab}\beta$ is the stabilizer of $\beta$ in $H$. Therefore, the $(t,m)$-placid
condition is equivalent to
\begin{equation}\lab{WB1}
2n(\beta)-\dim H/\mathrm{Stab}\beta > \frac12 (\dim W-2\dim H)
\end{equation}
for each $\beta\in\mathcal{B}$. In particular, this condition is satisfied when
$$2n(\beta)\ge \frac12\dim W$$ for all $\beta$. For example, if the set of weights is
symmetric with respect to the origin, the above is satisfied . This is the case for the
moduli spaces of vector bundles over a Riemann surface. (See Proposition \ref{modspex}.)

 When $K=U(1)$ acts on $M=\pp^n$ linearly and if $n_+, n_0, n_-$ denote the number of
positive, zero, negative weights respectively, then the condition (\ref{WB1}) is satisfied
if and only if $n_+=n_-$. Hence the $(t,m)$-placid condition may be viewed as a condition
on ``balancedness of weights''.

\begin{definition}\lab{def5.1k}\textrm{\cite{kw}}
The action on $M$ is said to be \emph{almost balanced}
if the condition (\ref{WB1}) is satisfied for all $\beta$ and $(H)$.
\end{definition}

\begin{remark}\lab{rem5.2c}
By (\ref{eq5.1c}), if almost balanced, we have an isomorphism
$$H^{<a_H}_H(W)\cong H^{<a_H}_H\left(W-\phi_W^{-1}(*)\right)$$
where $a_H=\frac12\dim W\git H=\frac12\dim W-\dim H$. This is an
easy consequence of the Gysin sequence (applied stratum by
stratum) because $\phi^{-1}_W(*)$ is stratified.
\end{remark}

\subsection{Almost-balanced action and GIT quotient}
We recall the following well-known facts from \cite{k1}: The
obvious action of $U(n+1)$ on $\pp^n$ is Hamiltonian with moment
map $\mu_{\pp^n}$. When $M\subset \pp^n$ is a smooth projective
variety and $K$ acts on $M$ via a homomorphism $K\to U(n+1)$, the
composition
$$\xymatrix{\mu:M\ar@{^(->}[r]& \pp^n\ar[r]^(.4){\mu_{\pp^n}}&
\mathfrak{u}(n+1)^*\ar[r]& \mathfrak{k}^*}$$ is the moment map for
$M$ and the set of semistable points in $M$ is equal to the
minimal Morse stratum $M^{ss}$ for $-|\mu|^2$. The GIT quotient
$M\git G$ of $M$ is homeomorphic to the symplectic reduction
$M\git K$ and the symplectic quotient map $\phi:M^{ss}\to M\git K$
is the GIT quotient map. Since we are interested in topology, we
will not distinguish symplectic quotients from GIT quotients.

\begin{theorem}\lab{th5.2}
Let $M\subset \pp^n$ be a smooth projective variety acted on by a
compact connected Lie group $K$ via a homomorphism $K\to U(n+1)$.
Suppose the $K$ action is almost balanced and there is at least
one point in $Z$ with finite stabilizer. Then we have a natural
map \begin{equation}\lab{eq5.2a}\phi_K ^*:IH^*(X)\to
H^*_K(M^{ss})\cong H^*_K(Z)\end{equation} of the middle perversity
intersection cohomology of $X=\mu^{-1}(0)/K$ into the
$K$-equivariant cohomology of $M^{ss}$. Moreover, $\phi^*_K$ is
injective and the intersection pairing of $IH^*(X)$ is given by
the cup product of $H^*_K(Z)$ from the formula
$$\phi_K ^*(\alpha)\cup \phi_K ^*(\beta)=\langle \alpha,\beta\rangle \phi_K ^*(\tau)$$
where $\tau$ is the class in $IH^{\dim X}(X)$ represented by a
point.
\end{theorem}
\begin{proof}
Certainly the GIT quotient map is subanalytic and the
stratification of $M^{ss}$ in \S\ref{subs4.2} is $K$-invariant.
The almost balanced condition is equivalent to the $(t,m)$-placid
condition and hence we have
$$\phi_K ^*:IH^*(X)\to H^*_K(M^{ss})\cong H^*_K(Z)$$
by Proposition \ref{p3.5}. For injectivity and intersection
pairing, it suffices to show that $\phi ^*_K(\tau)$ is nonzero in
$H^*_K(Z)$ by Corollary \ref{c3.6}.

Let $\Sigma$ denote the set of points in $Z$ whose stabilizer is
not finite. According to \cite{mes} \S4, there is an equivariant
proper map $\pi:\tilde{Z}\to Z$ such that
$\pi|_{\pi^{-1}(Z-\Sigma)}$ is a homeomorphism and the stabilizer
of every point in $\tilde{Z}$ is a finite group. Thus
$H^*_K(\tilde{Z})\cong H^*(\tilde{Z}/K)$.

Let $X^s=(Z- \Sigma)/K$ and consider the commutative diagram of natural maps
$$\begin{CD}
H^{\dim X}_c(X^s) @>{\phi ^*}>> H^{\dim X}_K(Z)\\
@V{\cong}VV         @VVV\\
H^{\dim X}(\tilde{Z}/K) @>{\cong}>> H^{\dim X}_K(\tilde{Z})
\end{CD}$$
where the subscript $c$ denotes compact support. The class $\tau$
is the image of a nonzero class in $H^{\dim X}_c(X^s)$ and hence
it follows from the above diagram that $\phi^*_K(\tau)$ is
nonzero.
\end{proof}

\begin{remark} Recall that (\ref{eq5.2a}) comes from a
morphism (\ref{eq3.4a}) in the derived category $\bold{D}^+_c(X)$
of bounded below cohomologically constructible sheaves over $X$
 because $\phi:M^{ss}\to X$
is $(t,m)$-placid. In particular, for any open subset $U\subset X$
we get a map $\phi^*_K:IH^*(U)\to H^*_K(\phi^{-1}(U))$ and it is
functorial with respect to restrictions.
\end{remark}


\section{The Kirwan map}
In this section, we recall the definition of the Kirwan map from
\cite{k3} and show that it is a left inverse of the pull-back
homomorphism $\phi_K^*$.

Let $M\subset \Bbb P^n$ be a connected nonsingular projective
variety acted on linearly by a connected reductive group $G$ via a
homomorphism $G\rightarrow GL(n+1)$. We may assume that the
maximal compact subgroup $K$ of $G=K^{\cc}$ acts unitarily
possibly after conjugation. Let $\mu:\Bbb P^n\rightarrow
u(n+1)^*\rightarrow \lik^*$ be the moment map for the action of
$K$. Then $M^{ss}=M\cap (\Bbb P^n)^{ss}$ retracts onto
$Z:=\mu^{-1}(0)\cap M$ by the gradient flow of $-|\mu|^2$ and the
GIT quotient $M\git G$ is homeomorphic to the symplectic quotient
$M\git K=M\cap \mu^{-1}(0)/K=:X$.

\subsection{Definition of the Kirwan map}\lab{subs6.1}
In order to define the Kirwan map, we assume that there is at
least one stable point in $M$, which amounts to saying that there
is at least one point in $Z$ whose stabilizer is finite.

We quote the following definitions from \cite{k2}.
\begin{definition}\begin{enumerate}
\item Let $\mathcal{R} (M)$ be a set of representatives of the conjugacy classes
of identity components of all subgroups of $K$ which appear
as the stabilizer of some point $x\in Z=\mu^{-1}(0)$.
\item Let $M^{ss}_H$ denote the set of those $x\in M^{ss}$ fixed by
$H\in {\mathcal{R}} (M)$.
\item Let $r(M)=\mathrm{max}\{\dim H|\, H\in \mathcal{R} (M)\}$.\end{enumerate}
\end{definition}

Of course, $M^{ss}_H$ is a smooth complex manifold.

The definition of the Kirwan map is by induction on $r(M)$. When
$r(M)=0$, the action of $K$ is locally free and $X$ is an
orbifold. Thus $\mathcal{IC}^{\cdot}_X\cong \cc_X$. The pull-back
morphism $\phi_K^*$ in this case is equal to the adjunction
morphism $\cc_X\to R{\phi_K}_*\cc$ which is an isomorphism by
\cite{BL} Theorem 9.1 (ii). The Kirwan map is defined as the
hypercohomology
$$\kappa_M^{ss}:H^*_K(M^{ss})\to IH^*(X)$$ of the inverse
$R{\phi_K}_*\cc\to \cc_X\cong \mathcal{IC}^{\cdot}_X $ of
$\phi_K^*$.

Now suppose $r(M)>0$. Let $\pi:\hatm\to M^{ss}$ be the blow-up of
$M^{ss}$ along the submanifold $$\bigsqcup_{\dim
H=r(M)}GM^{ss}_H.$$ If we choose a suitable linearization $\hatm$,
the semistable points in the closure of $\hatm$ all lie in $\hatm$
and $r(\hatm)<r(M)$. Let $\hat{X}=\hatm\git G$ and
$\hat{\phi}:\hatm^{ss}\to \hat{X}$ be the GIT quotient map. If we
keep blowing up in this fashion, we get a quasi-projective variety
$\widetilde{M}^{ss}$ birational to $M$ whose quotient
$\widetilde{X}$ has only finite quotient singularities. This is
called the \emph{partial desingularization} of $X$. See \cite{k2}
for details.

We have the following commutative diagram
\begin{equation}\lab{dia6.1}
\xymatrix{
{\hatm^{ss}}\ar@{^(->}[r]^{\imath}\ar[d]_{\hat{\phi}}&{\hatm}\ar[r]^{\pi}
&M^{ss}\ar[d]^{\phi}\\
{\hat{X}}\ar[rr]^{\sigma}& & X . }\end{equation} Since the GIT
quotient $\hat{X}$ is the categorical quotient of $\hatm^{ss}$,
$\sigma$ is defined uniquely by the universal property of the
categorical quotient. (Recall that GIT quotients are categorical
quotients \cite{mfk}.)

Inductively, we may suppose that we have a morphism
$$\kappa_{\hatm}^{ss}:R (\hat{\phi}_K) _* \cc \to
\mathcal{IC}^{\cdot}_{\hat{X}}.$$ Then we have a morphism
\begin{equation}\lab{eq6.2a} R{\phi_K}_*R{\pi_K}_*\cc\to
R{\phi_K}_*R{\pi_K}_*R{\imath_K}_*\cc=R\sigma_*R(\hat{\phi}_K)_*\cc\to
R\sigma_*\mathcal{IC}^{\cdot}_{\hat{X}}\end{equation} by composing
the above with the adjunction morphism $\cc\to
R{\imath_K}_*\imath_K^*\cc=R{\imath_K}_*\cc$, where
$\pi_K:EK\times_K\hatm\to EK\times_KM^{ss}$ is the induced map
from $\pi$ and $\imath_K$ is defined similarly. This induces a
homomorphism
\begin{equation}\lab{eq6.3a}\kappa^{ss}_{\hatm}\circ
{\imath_K}^*:H^*_K(\hatm)\to H^*_K(\hatm^{ss})\rightarrow
IH^*(\hatm\git G).\end{equation}

Next, compose (\ref{eq6.2a}) with the adjunction morphism
$$R{\phi_K}_*\cc\to
R{\phi_K}_*R{\pi_K}_*\pi_K^*\cc=R{\phi_K}_*R{\pi_K}_*\cc$$ to get
a morphism
\begin{equation}\lab{eqkeq}R {\phi_K} _* \cc \to R\sigma_*
\mathcal{IC}^{\cdot}_{\hat{X}}.\end{equation}

By \cite{k2} \S3, $\hat{X}$ is just the blow-up of $X$ along
$$\bigsqcup_{\dim H=r(M)}GM^{ss}_H\git G.$$ In particular, $\sigma$ is
proper and hence we can apply the decomposition theorem of
Beilinson, Bernstein, Deligne and Gabber in \cite{BBD} which says
\begin{equation}\lab{bbd6.3}R\sigma_*\mathcal{IC}^{\cdot}_{\hat{X}}=\mathcal{IC}^{\cdot}_X
\oplus \mathcal{F}^{\cdot}\end{equation}
 where
$\mathcal{F}^{\cdot}$ is a sheaf complex supported on the blow-up
center. Therefore, we have a morphism
\begin{equation}\lab{eqkic} R\sigma_*\mathcal{IC}^{\cdot}_{\hat{X}}\to
\mathcal{IC}^{\cdot}_X\end{equation} whose kernel is
$\mathcal{F}^{\cdot}.$

The composition of (\ref{eqkeq}) with (\ref{eqkic}) is the desired
morphism
$$\kappa_M^{ss}:R {\phi_K} _* \cc \to\mathcal{IC}^{\cdot}_X$$
which induces a homomorphism
$$H^*_K(\phi^{-1}(U))\to IH^*(U)$$
for an open set $U$ in $X$. This is denoted also by
$\kappa_M^{ss}$ by abuse of notations and called the \emph{Kirwan
map}.

\subsection{The pull-back is a right inverse}
 From \S \ref{almost}, when the
$K$ action on $M$ is almost balanced we have a morphism
$$\phi_K^*:\mathcal{IC}^{\cdot}_X\to R{\phi_K}_*\cc.$$
This induces a homomorphism $$IH^*(U)\to H^*_K(\phi^{-1}(U))$$ for
any open set $U$ in $X$ which we also denote by $\phi_K^*$ by
abuse of notations. The Kirwan map is a left inverse of
$\phi_K^*$.
\begin{proposition}\lab{prop6.2b} $\kappa_M^{ss}\circ \phi_K^*=1$.\end{proposition}
\begin{proof}
We know both $\phi_K^*$ and $\kappa_M^{ss}$ come from morphisms in
the derived category $\bold{D}_c^+(X)$. If we compose them, we get
a morphism
$$\mathcal{IC}^{\cdot}_X\to R {\phi_K} _* \cc \to \mathcal{IC}^{\cdot}_X.$$

On the set of stable points $M^s$ in $M^{ss}$, the action of $K$
is locally free (i.e. the stabilizers are finite groups). Let
$X^s=\phi(M^s)$ which is an orbifold. Then $\phi_K^*|_{X^s}$ is
the adjunction morphism
$$\mathcal{IC}^{\cdot}_{X^s}\cong \cc_{X^s}
\to R {\phi_K} _* \cc|_{X^s}$$ which is an isomorphism by
\cite{BL} Theorem 9.1 (ii) again, and $\kappa_M^{ss}|_{X^s}$ is
its inverse by definition since $X^s$ is untouched by the blow-ups
in the partial desingularization process. Therefore,
$\kappa_M^{ss}|_{X^s}\circ \phi_K^*|_{X^s}$ is the identity.

It is well-known (\cite{bo}, V \S9) that a morphism
$\mathcal{IC}^{\cdot}_X\to \mathcal{IC}^{\cdot}_X$ which restricts
to the identity over the smooth part (that is obviously contained
in $X^s$) is unique. Therefore, $\kappa_M^{ss}\circ \phi_K^*=1$.
\end{proof}

In particular, $\phi_K^*$ is injective and $\kappa_M^{ss}$ is
surjective.


\section{The image of $\phi_K^*$}\lab{section7}

In this section,  we identify the image of $\phi_K^*$  with a
naturally defined subspace $V^*_M\subset H^*_K(M^{ss})$ under a
slightly stronger assumption than almost balanced condition.

\subsection{Weakly-balanced action}
Let us make precise our assumption. We use the notations of
\S\ref{subs6.1}.

\begin{definition}\lab{firstWB}
Let $M\subset \pp^n$ be a projective variety with an action of a
compact Lie group $K$ via a homomorphism $K\to U(n+1)$. We say the
$K$ action on $M$ is \emph{weakly balanced} if it is almost
balanced and so is the $N^H/H$ action on the $H$-fixed submanifold
$M_H$ for each $H\in {\mathcal{R}} (M)$, where $N^H$ is the
normalizer of $H$ in $K$.
\end{definition}

For practical application, the following 2-step equivalent
definition is more useful. Recall that $G$ is the complexification
of $K$ which acts on $M$ via a homomorphism $G\to GL(n+1)$.

\begin{definition}\lab{5.3} \begin{enumerate}
\item Suppose a nontrivial compact group
$H$ acts on a vector space $W$ unitarily. Using the notations of
\S \ref{almost}, the action is said to be \emph{weakly linearly
balanced} if
 \begin{equation}\lab{ineq7.2}
 2n(\beta)- \dim  H/\mathrm{Stab}\beta > \frac12 (\dim W- 2\dim
 H)\end{equation}
for every $\beta \in \mathcal B$.

\item The $K$-action on $M$ is said to be \emph{weakly balanced}
 if for each $H\in \mathcal R(M)$ and for a point $x\in \mu^{-1}(0)$
 with $\mathrm{Lie}\,\mathrm{Stab}(x)=\mathrm{Lie} H$,
 the linear action of $H$ on the normal space $\mathcal N_x$
to $GM^{ss}_H$ is weakly linearly balanced and so is the action of
$(H\cap N^{L})/L$ on the $L$-fixed linear subspace $\mathcal
N_{x,L}$ for each connected subgroup $L$ of $H$ whose conjugate
appears in $\mathcal{R}(M)$.
\end{enumerate}
\end{definition}

\begin{lemma}\lab{lem7.3a}
The two definitions \ref{firstWB} and \ref{5.3} are
equivalent.\end{lemma}

\begin{proof}

Let $x\in \mu^{-1}(0)$ and $P=\mathrm{Stab}(x)$. Let $L$ be the
identity component of $P$. By Lemma \ref{local}, a neighborhood of
$Kx$ is equivariantly diffeomorphic to
$$K\times_P\left((\lik/\mathfrak{p})^*\times W\right)\times \rr^{\dim
X_{(P)}}$$ for some symplectic $P$-vector space $W$. Let us call
$W$ the \emph{normal slice} at $x$.

Since $P/L$ is discrete, $P\subset N^L$ and $P$ acts on the
$L$-fixed subspace $W^L$ of $W$. By direct computation, one can
check that $M^{ss}_L$ in this neighborhood is
$$N^L\times_P\left((\mathfrak{n}^L/\mathfrak{p})^*\times
W^L\right)\times \rr^{\dim X_{(P)}}$$ and hence $GM^{ss}_L$ is
$$K\times_P\left((\lik/\mathfrak{p})^*\times W^L\right)\times \rr^{\dim
X_{(P)}}$$ where $\mathfrak{n}^L$ (resp. $\mathfrak{p}$) is the
Lie algebra of $N^L$ (resp. $P$). Therefore the normal space
$\cN_x$ to $GM^{ss}_L$ is the orthogonal complement of $W^L$ in
$W$.

 If $x$ is a generic point in $M^{ss}_L\cap
\mu^{-1}(0)$ such that $P$ is minimal among those containing $L$,
then $\mathcal N_x=W$. Suppose the $K$-action on $M$ is almost
balanced, i.e. the action of $\mathrm{Stab}(x)$ on the normal
slice $W$ at $x$ is weakly linearly balanced for all $x\in
\mu^{-1}(0)$. Then by choosing a generic $x$ for each $L\in
\mathcal{R}(M)$, we deduce that the action of $L$ on the normal
space $\cN_x=W$ to $GM^{ss}_L$ is weakly linearly balanced.

In general, we only have $\mathcal N_x\subset W$. But $L$ acts
trivially on $\mathcal N_x^{\perp}\cap W$ and hence the weights of
the maximal torus action on $\mathcal N_x^{\perp}\cap W$ are all
zero. By examining the inequality (\ref{ineq7.2}) it is easy to
see that if the $L$-action on $\mathcal N_x$ is weakly linearly
balanced then so is the $P$ action on $W$. Therefore, if for each
$L\in \mathcal{R}(M)$ the action of $L$ on the normal space to
$GM^{ss}_L$ at a generic point $x\in \mu^{-1}(0)$ with
$\mathrm{Lie}\mathrm{Stab}(x)=\mathrm{Lie}(L)$ is weakly linearly
balanced, then the action of $\mathrm{Stab}(x)$ on the normal
slice $W$ at $x$ is weakly linearly balanced for all $x\in
\mu^{-1}(0)$, i.e. the $K$-action on $M$ is almost balanced.

Now let $J\in\mathcal{R}(M)$ and suppose $x\in M^{ss}_J\cap
\mu^{-1}(0)$. Then $L\supset J$. Using Lemma \ref{lem8.1a} which
is purely a group theoretic result, it is direct to check that in
the neighborhood of $Kx$, $M^{ss}_J$ is
$$\begin{array}{lc} &N^{J}\times_{P\cap
N^{J}}\left(
(\mathfrak{n}^{J}/\mathfrak{p}\cap\mathfrak{n}^{J})^*\times W^{J}
\right)\times\rr^{\dim X_{(P)}}\\ & \cong N^{J}/J\times_{P\cap
N^{J}/J}\left(
(\mathfrak{n}^{J}/\mathfrak{p}\cap\mathfrak{n}^{J})^*\times W^{J}
\right)\times\rr^{\dim X_{(P)}}\end{array}
$$
See the proof of Proposition \ref{prop8.3a} for a similar
computation.

If $\mathrm{Stab}(x)$ is minimal with
$\mathrm{Lie}\mathrm{Stab}(x)=\mathrm{Lie}L$ so that $\cN_x=W$,
then the $J$-fixed sets $\cN_{x,J}$ and $W^J$ are isomorphic. In
general, we only have $\cN_{x,J}\subset W^J$. But by the arguments
in the previous paragraphs, we deduce that the action of $N^J/J$
on $M_J$ is almost balanced if and only if the action of $L\cap
N^J/J$ on $\cN_{x,J}$ is weakly linearly balanced for $x\in
\mu^{-1}(0)$ with $\mathrm{Lie}\mathrm{Stab}(x)=\mathrm{Lie}L$ for
all $L\supset J$. So we proved the lemma.
\end{proof}

The weakly balanced condition is satisfied by many interesting
spaces including the diagonal $SL(2)$ action on $(\pp^1)^{2n}$.
(See \S9.) Also, it is satisfied by (the GIT construction of) the
moduli spaces of holomorphic vector bundles over a Riemann surface
of any rank and any degree. For the next proposition, let us use
Definition \ref{5.3}.

\begin{proposition} \lab{modspex}
Let $M(n,d)$ be the moduli space of rank $n$ holomorphic vector
bundles of degree $d>n(2g-1)$ over a Riemann surface $\Sigma$ of
genus $g$, which is a GIT quotient of a nonsingular
quasiprojective variety $\frak R(n,d)^{ss}$ by $G=SL(p)$ for
$p=d-n(g-1)$. (See \cite{New, k5}.) The action of $SL(p)$ on
$\frak R(n,d)^{ss}$ is weakly balanced.\end{proposition}
\begin{proof} Let $E$ be a semistable vector bundle
such that $E\cong m_1E_1\oplus \cdots \oplus m_sE_s$ where $E_i$'s
are non-isomorphic stable bundles with the same slope. Then the
identity component of $\mathrm{Stab} \,E$ in $G$ is
$H^{\cc}=S(\prod_{i=1}^sGL(m_i))$ where $S$ denotes the subset of
elements whose determinant is 1. The normal space to $GM^{ss}_H$
at $E$ is (\cite{ab2,k5})
$$\begin{array}{lll} H^1(\Sigma, End'_{\oplus}E)&=H^1(\Sigma, \oplus_{i,j}
(m_im_j-\delta_{ij})Hom(E_i,E_j))\\
&=\oplus_{i,j}H^1(\Sigma, (m_im_j-\delta_{ij})E^*_i\otimes E_j)\end{array}$$
More precisely,
$$\begin{array}{lll}
 H^1(\Sigma, End'_{\oplus}E)=&\oplus_{i<j}[H^1(\Sigma,E^*_i\otimes E_j)
\otimes Hom(\Bbb C^{m_i}, \Bbb C^{m_j}) \\
&\oplus H^1(\Sigma,E_i\otimes E^*_j)
\otimes Hom(\Bbb C^{m_j}, \Bbb C^{m_i})]\\
&\oplus [\oplus_iH^1(\Sigma, End\,E_i)\otimes sl(m_i)]\end{array}$$
Because $E_i$ is not isomorphic to $E_j$ for $i\ne j$, $H^0(\Sigma,
E^*_i\otimes E_j)=0=H^0(\Sigma,E_i\otimes E^*_j)$ and thus
$$\begin{array}{lll}
 \dim H^1(\Sigma, E^*_i\otimes E_j)&=-RR(\Sigma, E^*_i\otimes E_j)
=(rank\,E_i)(rank\, E_j)(g-1)\\
&=-RR(\Sigma, E_i\otimes E^*_j)=\dim H^1(\Sigma, E_i\otimes
E^*_j)\end{array}$$ where $RR$ denotes the Riemann-Roch number.
Therefore, the weights of the representation of $H^{\cc}$ on
$H^1(\Sigma, End'_{\oplus}E)$ are symmetric with respect to the
origin. This implies that the action is weakly linearly balanced.
As each subgroup $L^{\cc}$ as in Definition \ref{5.3} (2) is
conjugate to $S(\prod_{i=1}^sGL(m'_i))$ for a ``subdivision''
$(m_1',m_2',...)$ of $(m_1,m_2,...)$, it is easy to check that
such $H\cap N^{L}/L$ action on the $L$-fixed point set is also
weakly linearly balanced. \end{proof}

\subsection{The image of the pull-back homomorphism}

For any $H\in \mathcal{R}(M)$, consider the natural map (sometimes
called the ``resolution'')
\begin{equation}
K\times _{N^H} M^{ss}_H\rightarrow KM^{ss}_H\end{equation} and the
corresponding map on the cohomology (\cite{k5} Lemma 1.21)
\begin{equation}\lab{D3.1}
H^*_K(KM^{ss}_H)\rightarrow H^*_K(K\times _{N^H} M^{ss}_H)
\cong H^*_{N^H}(M^{ss}_H)
\cong [H^*_{N^H_0/H}(M^{ss}_H)\otimes H^*_H]^{\pi_0N^H}
\end{equation}
where $N_0^H$ is the identity component of $N^H$. For any
$\zeta\in H^*_K(M^{ss})$ let $\zeta|_{K\times _{N^H} M^{ss}_H}$
denote the image of $\zeta$ by the composition of the above map
with the restriction map $H^*_K(M^{ss})\rightarrow
H^*_K(KM^{ss}_H)$. Now, we can describe the image of $\phi_K^*$.

\begin{definition}\lab{def7v}
Put $n_H=\frac12 \codim GM^{ss}_H -\dim H$ and
$H^{<n_H}_H=\oplus_{i<n_H}H^i_H$. We define $V^*_M$ as the set of
$\zeta\in H^*_K(M^{ss})$ such that
\begin{equation}\lab{defVtensor} \zeta|_{K\times_{N^H}M^{ss}_H}
\in H^*_{N^H_0/H}(M_H^{ss}) \otimes H^{<n_H}_H\end{equation} for
each $H\in \mathcal{R}(M)$.
\end{definition}

\begin{remark}
The definition of $V_M^*$ is independent of the choices of $H$s in
the conjugacy classes and the tensor product expressions: The
former is easy to check by translating by $g$ if $H$ is replaced
by $gHg^{-1}$. The latter can be immediately seen by considering
the gradation of the degenerating spectral sequence for the
cohomology of the fibration
\begin{equation}\lab{D3.3}
\begin{CD}
(EN_0^H\times EN^H_0/H)\times_{N^H_0}M^{ss}_H\\
      @VVV\\
        EN^H_0/H\times_{N^H_0/H}M^{ss}_H\end{CD}
\end{equation}
The fiber is homotopically equivalent to $BH$. (See \cite{k5}
Lemma1.21.) Though the last isomorphism in (\ref{D3.1}) is not
canonical, the subspace in (\ref{defVtensor}) is
canonical.\end{remark}

From (\ref{defVtensor}), we have
\begin{equation}\lab{defVker}
V^*_M=\mathrm{Ker}\left( H^*_K(M^{ss})\to \bigoplus_{H\in
\mathcal{R}(M)}H^*_{N^H_0/H}(M^{ss}_H)\otimes H^{\ge n_H}_H\right)
\end{equation} and thus $V_M^*$ can be thought of as a subset of
$H_K^*(M^{ss})\cong H^*_K(Z)$, obtained by ``truncating locally''.

Now we can state the main theorem of the section which will be proved
in the next section.
\begin{theorem}\lab{splittingtheorem}
Let $M\subset\pp^n$ be a projective smooth variety acted on
unitarily by a compact connected group $K$ with at least one
stable point. Suppose that the weakly balanced condition is
satisfied. Then we have $\phi_K^*(IH^*(X))=V^*_M$. Moreover, for
any open set $U$ of $X$, if we define $V^*_{\phi^{-1}(U)}\subset
H^*_K(\phi^{-1}(U))$ as in Definition \ref{def7v}, then we have
$\phi_K^*(IH^*(U))=V^*_{\phi^{-1}(U)}$.
\end{theorem}


\section{Proof of Theorem \ref{splittingtheorem}}
This section is devoted to a proof of Theorem
\ref{splittingtheorem}. Let us use the notations of
\S\ref{subs6.1} and \S\ref{section7}. Recall that
\begin{equation}\lab{D4.1}
r:=r(M)=\mathrm{max}\{\dim H\,|\, H\in \mathcal
R(M)\}.\end{equation} Our proof is by induction on $r(M)$.

When $r=0$, we have nothing to prove since
$$V^*_M=H^*_K(M^{ss})\cong IH^*(M\git G).$$
So we consider the case $r>0$. Suppose the theorem is true for all
projective varieties $\Gamma$ with $r(\Gamma)\le r-1$. Let $\hatm$
be the blowup of $M^{ss}$ along the submanifold $$\bigsqcup_{\dim
H=r}GM^{ss}_H.$$ Then from \cite{k2} \S6, we have
$${\mathcal{R}}(\hatm)=\{H\in \mathcal R(M)\, |\, \dim H\le
r-1\}$$ and thus $r(\hatm)\le r-1$.

For simplicity, we assume from now on that there exists only one
$H$ such that $\dim H=r$. (The general case is no more difficult
except for repetition. We can deal with each $H$ one by one. See
\cite{k2}, Cor.8.3.) We fix this $H$ once and for all till the end
of this section.

\begin{remark}\lab{remWBhat}
(1) To be precise, we have to take the closure of $\hatm$ with
respect to a suitable linearization described in \cite{k2} and
then resolve the possible singularities. But as argued in
\cite{k2}, this does not cause any trouble for us because all the
semistable points are contained in $\hatm$ and we are only
interested in the semistable points.

(2) By \cite{k3} 1.6, for $L\in \mathcal{R}(\hatm)$, the $L$-fixed
set $\hatm^{ss}_L$ in $\hatm^{ss}$ is the proper transform of the
$L$-fixed set $M^{ss}_L$ in $M^{ss}$. In particular, the normal
space to $G\hatm^{ss}_L$ in $\hatm^{ss}$ at a generic point is
isomorphic to the normal space to $GM^{ss}_L$ in $M^{ss}$ at a
generic point. Notice that the weakly balanced condition in
Definition \ref{5.3} is purely about the actions of $L$ on the
normal spaces $\cN_x$ to $GM^{ss}_L$ for $L\in\mathcal{R}(M)$.
(The fixed set by a subgroup of $L$ is determined by the action of
$L$.) Therefore, if the $K$-action on $M^{ss}$ is weakly balanced,
the action on $\hatm$ is also weakly balanced.
\end{remark}

By our induction hypothesis, $\hat{\phi}_K^*(IH^*(\hatm\git
G))=V^*_{\hatm}$ and the same holds for any open subset of
$\hat{X}=\hatm\git G$.

We start the proof with a few lemmas.

\begin{lemma}\lab{lem8.1a}(\cite{k2} Proposition 8.10) Suppose
$L\subset P$ are compact subgroups of $K$ and $L$ is connected.
Then there exist finitely many elements $k_1,k_2,\cdots,k_m$ in
$K$ such that $$\{k\in K\,|\,k^{-1}Lk\subset P\}=\bigsqcup_{1\le
i\le m}N^Lk_iP$$ where $N^L$ is the normalizer of $L$ in
$K$.\end{lemma}
\begin{proof}
See the proof of \cite{k2} p77.\end{proof}

Let $E$ be the exceptional divisor in $\hatm$ of the blow-up.
\begin{lemma}\lab{lem8.2a} By restriction,
we have an isomorphism $$\mathrm{Ker}\left(H^*_K(\hatm)\to
H^*_K(\hatm^{ss})\right)\cong \mathrm{Ker}\left( H^*_K(E)\to
H^*_K(E^{ss}) \right).$$
\end{lemma}
\begin{proof} This follows from \cite{k2} 7.5, 7.6 and 7.11.\end{proof}

We first show that the image of $\phi_K^*$ is contained in
$V^*_{M}$. Let $U$ be an open subset of $X=M\git G$.
\begin{proposition} \lab{prop8.3a}
$\phi^*_K(IH^*(U))\subset V^*_{\phi^{-1}(U)}$.\end{proposition}
\begin{proof}
Recall that $E\mathcal{G}$ denotes a contractible free
$\mathcal{G}$-space for a Lie group $\mathcal{G}$ and
$B\mathcal{G}=E\mathcal{G}/\mathcal{G}$. Let $L\in \mathcal R(M)$.
From the obvious commutative diagram
$$\xymatrix{
\big(EK\times
E(N^{L}/L)\big)\times_{N^{L}}M^{ss}_L\ar[r]^(.65){g}\ar[d]_{\phi'_L}
&EK\times_KM^{ss}\ar[d]^{\phi_K}\\
E(N^{L}/L)\times_{N^{L}/L}M^{ss}_L\ar[r]^(.7){f}&X }$$ we get a
morphism
$$f^*R{\phi_K}_*\cc\to R{\phi'_L}_*g^*\cc\to \tau_{\ge
n_L}R{\phi'_L}_*g^*\cc=\tau_{\ge n_L}R{\phi'_L}_*\cc.$$ This
induces a morphism by adjunction
\begin{equation}\lab{rhoeq8.3} R{\phi_K}_*\cc\to Rf_*f^*R{\phi_K}_*\cc\to
Rf_*\tau_{\ge
n_L}R{\phi'_L}_*\cc=:\mathcal{A}^{\cdot}_L.\end{equation} The
fiber of $\phi'_L$ is homotopically equivalent to $BL$ and thus
this morphism induces the truncation homomorphism
$$H^*_K(M^{ss})\to [H^*_{N^L_0/L}(M^{ss}_L)\otimes H^{\ge n_L}_L
]^{\pi_0N^L}$$ where $N^L_0$ is the identity component of $N^L$.
By composing (\ref{rhoeq8.3}) with $\phi_K^*$, we get a morphism
$$\xymatrix{\rho:\mathcal{IC}^{\cdot}_X\ar[r]^{\phi_K^*}&R{\phi_K}_*\cc\ar[r]&
\mathcal{A}^{\cdot}_L }$$ whose hypercohomology gives us
$$IH^*(U)\to H^*_K(\phi^{-1}(U))\to  [H^*_{N^L_0/L}(M^{ss}_L\cap
 \phi^{-1}(U))\otimes H^{\ge n_L}_L
]^{\pi_0N^L}.$$ Therefore it suffices to show that $\rho$ is equal
to zero in view of (\ref{defVker}).

The sheaf complex $\mathcal{A}^{\cdot}_L$ is trivial on the
complement of the closed subset $GM^{ss}_L\git G$. Hence $\rho$ is
zero on this open dense subset. Hence by adding stratum by stratum
in the order of increasing codimension, it suffices to show the
following: Let $P$ be a subgroup of $K$ and consider the stratum
$X_{(P)}$ defined in \S\ref{section4}. Suppose $U$ is an open
subset of $X$ containing $X_{(P)}$ such that $U-X_{(P)}$ is open
and  $\rho|_{U-X_{(P)}}$ is equal to zero. Then $\rho|_U$ is also
zero.

Let $\imath:U-X_{(P)}\hookrightarrow U$ and put $n_P=\frac12\codim
X_{(P)}$. We claim that
\begin{equation}\lab{eqcl8.3}
\tau_{<n_P}\mathcal{A}^{\cdot}_L|_U\cong
\tau_{<n_P}R\imath_*\mathcal{A}^{\cdot}_L|_{U-X_{(P)}}.\end{equation}
This claim enables us to deduce that $\rho|_U$ is zero from
$\rho|_{U-X_{(P)}}$ being zero because $\rho|_U$ is the
composition
$$\xymatrix{
{\mathcal{IC}_X^{\cdot}|_U}\ar[r]^(0.3){\cong}&
{\tau_{<n_P}R\imath_*\mathcal{IC}_X^{\cdot}|_{U-X_{(P)}}}
\ar[r]^0&{\tau_{<n_P}R\imath_*\mathcal{A}_L^{\cdot}|_{U-X_{(P)}}}\ar[r]^(0.6){\cong}&
{\tau_{<n_P}\mathcal{A}^{\cdot}_L}|_U\ar[r]&{\mathcal{A}^{\cdot}_L}|_U}$$
which is zero.

Let us now prove (\ref{eqcl8.3}). If $L$ is not conjugate to a
subgroup of $P$, then $X_{(P)}$ does not intersect with
$GM^{ss}_L\git G$ and thus we have nothing to prove. So we may
assume $L\subset P$ after conjugation if necessary.

Consider the commutative diagram
$$\xymatrix{
M^{ss}_L\ar@{^(->}[r]^{\jmath}\ar[d]^{q}&M^{ss}\ar[d]^{\phi}\\
M^{ss}_L\git N^L\ar[r]^h&M\git K=X }$$ where $h$ is the unique map
defined by the universal property of the categorical quotient
$M^{ss}_L\git N^L$ of $M^{ss}_L$.

We compute the stalk cohomology of both sides of (\ref{eqcl8.3}).
By Lemma \ref{local}, the preimage of a contractible neighborhood
$\Delta$ of a point in $X_{(P)}$ by $\phi$ is equivariantly
homeomorphic to
\begin{equation}\lab{eq8.3a} K\times_P\left(
(\mathfrak{k}/\mathfrak{p})^*\times W\right)\times \rr^{\dim
X_{(P)}}\end{equation} for some symplectic $P$-vector space $W$.
By Lemma \ref{lem8.1a}, it is direct to check that $M^{ss}_L$ in
this neighborhood is
\begin{equation}\lab{eqp8.3.5}\bigsqcup_{1\le i\le m}k_iN^{L_i}\times_{P\cap
N^{L_i}}\left(
(\mathfrak{n}^{L_i}/\mathfrak{p}\cap\mathfrak{n}^{L_i})^*\times
W^{L_i} \right)\times\rr^{\dim X_{(P)}}\end{equation} where
$L_i=k_i^{-1}Lk_i$ and $W^{L_i}$ is the $L_i$-fixed subspace of
$W$. Also $\mathfrak{n}^L$ denotes the Lie algebra of $N^L$.

If we delete $X_{(P)}$ from the neighborhood $\Delta$, then the
preimage by $\phi$ is
$$K\times_P\left( (\mathfrak{k}/\mathfrak{p})^*\times
(W-\phi_W^{-1}(*))\right)\times \rr^{\dim X_{(P)}}$$ where
$\phi_W:W\to W\git P$ is the GIT quotient map and $*$ is the
vertex of the cone $W\git P$. The intersection of this with
$M^{ss}_L$ is homeomorphic to
$$\bigsqcup_{1\le i\le m}k_iN^{L_i}\times_{P\cap N^{L_i}}\left(
(\mathfrak{n}^{L_i}/\mathfrak{p}\cap\mathfrak{n}^{L_i})^*\times
(W^{L_i}-\phi_W^{-1}(*)) \right)\times\rr^{\dim X_{(P)}}.$$

Hence the stalk cohomology of the left hand side of
(\ref{eqcl8.3}) is
$$\bigoplus_i H^{<n_P-n_L}_{P\cap N^{L_i}/L_i}(W^{L_i})\otimes H^{\ge n_L}_{L_i}$$
while the right hand side has
$$\bigoplus_i H^{<n_P-n_L}_{P\cap N^{L_i}/L_i}(W^{L_i}-\phi_W^{-1}(*))\otimes H^{\ge n_L}_{L_i}.$$
Thus it suffices to show that
\begin{equation}\lab{eqp8.3.1}H^{<n_P-n_L}_{P\cap N^{L_i}/L_i}(W^{L_i})\cong
H^{<n_P-n_L}_{P\cap N^{L_i}/L_i}(W^{L_i}-\phi_W^{-1}(*)).
\end{equation}
Without loss of generality, we may assume $L_i=L$.

By definition, we have
\begin{equation}\lab{eqdim8.6.1}
n_P=\frac12\dim W\git P=\frac12 \dim W-\dim P.\end{equation}
From
(\ref{eqp8.3.5}), it is easy to deduce that $GM^{ss}_L$ in the
preimage of $\Delta$ is
\begin{equation}\lab{eq8.6} K\times_P\left( (\mathfrak{k}/\mathfrak{p})^*\times
P_{\cc}W^L\right)\times \rr^{\dim X_{(P)}}\end{equation} where
$P_{\cc}$ is the complexification of $P$ in $G$ and $W$ is
assigned a complex structure compatible with the symplectic
structure. Hence, we have
\begin{equation}\lab{eqdim8.6.2}
n_L=\frac12 \codim GM^{ss}_L-\dim L=\frac12 \left(\dim W-\dim
P_{\cc}W^L\right)-\dim L.\end{equation} From the surjectivity of
the morphism
$$P_{\cc}\times _{P_{\cc}\cap N^L_{\cc}}W^L\to P_{\cc}W^L$$
we see that
\begin{equation}\lab{eqdim8.6.3}
\dim P_{\cc}W^L\le \dim P_{\cc}+\dim W^L-\dim (P_{\cc}\cap
N^L_{\cc}).\end{equation}

Comparing (\ref{eqdim8.6.1}), (\ref{eqdim8.6.2}) and
(\ref{eqdim8.6.3}), we get
\begin{equation}\lab{eqdim8.6.4}
n_P-n_L\le \frac12\dim W^L-\dim (P\cap N^L/L).\end{equation} Since
the action of $N^L/L$ on $M^{ss}_L$ is almost balanced (Remark
\ref{rem5.2c}),  by (\ref{eqdim8.6.4}) we have
\begin{equation}\lab{eqp8.3.2} H^{<n_P-n_L}_{P\cap
N^{L}/L}(W^{L})\cong H^{<n_P-n_L}_{P\cap
N^{L}/L}(W^{L}-\phi_{W^L}^{-1}(*))
\end{equation}
where $\phi_{W^L}:W^L\to W^L\git P\cap N^L$ is the GIT quotient
map.

Finally, we observe that
\begin{equation}\lab{eqp8.3.3}\phi_{W^L}^{-1}(*)=\phi_W^{-1}(*)\cap W^L.
\end{equation} This is because
we know the following from \cite{k1}:
\begin{enumerate}
\item For $x\in W^L$, $x\in\phi_{W^L}^{-1}(*)\Leftrightarrow
\lim_{t\to \infty}x_t=0$ where $x_t$ is the gradient flow for
$-|\mu_{W^L}|^2$ with $x_0=x$ ($\mu_{W^L}$ is the moment map for
$W^L$).
\item For $x\in W$, $x\in\phi_{W}^{-1}(*)\Leftrightarrow
\lim_{t\to \infty}x_t=0$ where $x_t$ is the gradient flow for
$-|\mu_{W}|^2$ with $x_0=x$ ($\mu_W$ is the moment map for $W$).
\item For a moment map $\mu$ on a symplectic manifold,
the gradient vector at $x$ for $-|\mu|^2$ is $-2i\mu(x)_x$ if
$\lik$ is identified with $\lik^*$ by the Killing form.
\item $\mu_W(x)\in \mathfrak{n}^L$ if $x\in W^L$ and hence
$\mu_W(x)=\mu_{W^L}(x)$.
\end{enumerate}
(\ref{eqp8.3.1}) follows from (\ref{eqp8.3.2}) and
(\ref{eqp8.3.3}).

\end{proof}

We need a few more lemmas.

\begin{lemma} \lab{lem8.5a} Consider the diagram (\ref{dia6.1}) in
\S\ref{subs6.1}. The restriction to $V^*_M$ of
$$\xymatrix{H^*_K(M^{ss})\ar@{^(->}[r]^(.6){\pi^*_K} &H^*_K(\hatm)\ar[r]^{\imath_K^*}& H^*_K(\hatm^{ss})}$$
factors through $V^*_{\hatm}$ and is injective. A similar
statement is true for $\phi^{-1}(U)$ where $U$ is any open set in
$X$.
\end{lemma}
\begin{proof} Let $\zeta$ be a nonzero element in $V_M^*$. Then
$$\zeta|_{K\times _{N^{L}}M^{ss}_{L}}
\in [H^*_{N_0^{L}/L}(M^{ss}_{L})\otimes H^{<n_L}_{L}]
^{\pi_0N^{L}}$$  for each  $L\in \mathcal R(M).$ Its image in
$H^*_K(\hatm^{ss})$ satisfies
$$
\zeta|_{K\times_{N^{L}}\hatm^{ss}_{L}} \in
[H^*_{N_0^{L}/L}(\hatm^{ss}_{L})\otimes H^{<n_L}_{L}]
^{\pi_0N^{L}}$$  for each  $L\in \mathcal R(\hatm)=\{L\in \mathcal
R(M)| \dim L<r(M)\}.$ This follows from the commutative diagram
$$\begin{CD}
\hatm^{ss}@<<< K\hatm^{ss}_{L} @<<< K\times_{N^{L}}\hatm^{ss}_{L}\\
@VVV   @VVV  @VVV\\
M^{ss}@<<< KM^{ss}_{L} @<<< K\times _{N^{L}}M^{ss}_{L}.
\end{CD}$$
Therefore, $\zeta$ is mapped to an element in $V_{\hatm}^*$.

Recall that $H$ is the identity component of a stabilizer which
has the maximal dimension $r(M)$ and the blowup center is the
submanifold $GM^{ss}_H$. Since $\pi^*_K$ is an injection by the
well-known argument in \cite{GH} p605, we may think of $\zeta$ as
an element of $H^*_K(\hatm)$. By Lemma \ref{lem8.2a}, if
$\zeta|_{GM^{ss}_H}=0$ i.e. $\zeta|_E=0$, then
$\zeta|_{\hatm^{ss}}\ne 0.$

Let us now consider the case when $\zeta|_{GM^{ss}_H}\ne 0$. Since
$H$ is maximal, we have an isomorphism
$$H^*_{N_0^{H}/H}(M^{ss}_{H})\cong H^*(M_H^{ss}\git N^H_0).$$
By the definition of $V_M^*$, $\zeta|_{GM^{ss}_H}$ lies in
\begin{equation}\lab{isomNb}
[H^*(M_H^{ss}\git N^H_0)\otimes
H^{<n_H}_H]^{\pi_0N^H}\end{equation}
and we have
\begin{equation}\lab{isomNa}
H^*_K(E^{ss})=[H^*(M_H^{ss}\git N^H_0)\otimes H^*_H(\Bbb P\mathcal
N_x^{ss})]^{\pi_0N^H}\end{equation} from \cite{k3} Lemma 1.16
where $\cN_x$ is the normal space to $GM^{ss}_H$. Because the
$K$-action is almost balanced, the codimensions of the unstable
strata in $\pp\mathcal N_x$ are greater than $n_H$ by (\ref{WB1}).
Therefore, by the equivariant Morse theory \cite{k1}, we deduce
that the restriction homomorphism
$$H^{<n_H}_H(\pp\mathcal N_x)\to H^{<n_H}_H(\pp\mathcal N_x^{ss})$$
is an isomorphism. In particular, $H^{<n_H}_H$ injects into
$H^{<n_H}_H(\pp\mathcal N_x^{ss})$. By (\ref{isomNa}) and
(\ref{isomNb}), the image of $\zeta$ in $H^*_K(E^{ss})$ is not
zero and thus $\zeta$ injects into $H^*_K(\hatm^{ss})$.

It is obvious from our proof that the statement is true for any
open set $U$ in $X$.
\end{proof}

Let $\cN$ be the normal bundle to $GM^{ss}_H$  in $M^{ss}$ and
$\cN_U$ be the normal bundle to $GM^{ss}_H\cap \phi^{-1}(U)$ for
any open subset $U$ of $X$. It is proved in \cite{k3} Lemma 2.9
that $\cN\git K$ is homeomorphic to a neighborhood of
$X_{(\lih)}=GM^{ss}_H\git G$ and hence $\cN_U\git K$ is
homeomorphic to a neighborhood of $U\cap X_{(\lih)}$. We identify
$\cN_U$ with a tubular neighborhood of $GM^{ss}_H\cap
\phi^{-1}(U)$ and identify $\cN_U\git K$ with a neighborhood $U_1$
of $U\cap X_{(\lih)}$ in $U$.

By the gradient flow of $-|\mu|^2$, $M^{ss}$ can be equivariantly
retracted into $\mu^{-1}(D_{\varepsilon})$ where $D_{\varepsilon}$
is the disk of radius $\varepsilon$ around $0$ in $\lik^*$. By
shrinking $U_1$ if necessary and taking $\varepsilon$ sufficiently
small, $\phi^{-1}(U_1)$ is retracted into $\cN_U$. Conversely, if
we decrease the radius of a tubular neighborhood of $GM^{ss}_H\cap
\phi^{-1}(U)$, it is included in $\phi^{-1}(U_1)$. These two
inclusions are clearly inverse to each other \emph{homotopically}
and $K$-equivariantly since $\mu$ is equivariant. In particular,
$\cN_U$ and $\phi^{-1}(U_1)$ are homotopically equivalent open
neighborhoods of $GM^{ss}_H\cap \phi^{-1}(U)$. Therefore,
$H^*_K(\phi^{-1}(U_1))$ is canonically isomorphic to
$H^*_K(\cN_U)$ and $V^*_{\phi^{-1}(U_1)}\cong V^*_{\cN_U}$. Hence
for cohomological purpose, we can think of $\cN_U$ as the preimage
of a neighborhood of $U\cap X_{(\lih)}$.

\begin{proposition}\lab{lem8.6a} The restriction of the Kirwan map
gives us an isomorphism $V^*_{\cN_U}\cong IH^*(\cN_U\git G)$.
\end{proposition}

We postpone the proof of this proposition and prove Theorem
\ref{splittingtheorem}.

\begin{proof}[Proof of Theorem \ref{splittingtheorem}] For
an open subset $U$ of $X$, let $B^*(U)$ be the kernel of the
Kirwan map restricted to $V^*_{\phi^{-1}(U)}$. Then by Proposition
\ref{prop8.3a}, we can write
\begin{equation}\lab{vdec8}V^*_{\phi^{-1}(U)}=IH^*(U)\oplus B^*(U).\end{equation}
We have to show that $B^*(U)$ is zero.

Let $\hat{U}$ be the preimage of $U$ by the blow-up map
$\sigma:\hat{X}\to X$. By our induction hypothesis, the pull-back
$\hat{\phi}_K^*$ is an isomorphism of $IH^*(\hat{U})$ onto
$V^*_{\hat{\phi}^{-1}(\hat{U})}$ and the Kirwan map is its
inverse.

Recall that we have the decomposition (\ref{bbd6.3}) which induces
an isomorphism \begin{equation}\lab{eqtp8}IH^*(\hat{U})\cong
IH^*(U)\oplus F^*(U)\end{equation} where $F^*(U)$ is the
hypercohomology of $\mathcal{F}^{\cdot}$ over $U$. Since $B^*(U)$
is mapped to zero by
$$V^*_{\phi^{-1}(U)}\hookrightarrow V^*_{\hat{\phi}^{-1}(\hat{U})}\cong
IH^*(\hat{U})\cong IH^*(U)\oplus F^*(U)\to IH^*(U)$$ we see that
$B^*(U)$ injects into $F^*(U)$.

As $\phi_K^*$, $\kappa_M^{ss}$ and (\ref{eqtp8}) all came from
sheaf complexes, we have the following commutative diagram by
restriction
$$\xymatrix{
V^*_{\phi^{-1}(U)}\ar@{^(->}[r]\ar[d] &
V^*_{\hat{\phi}^{-1}(\hat{U})}\ar[r]^{\cong}\ar[d] &
IH^*(\hat{U})\ar[d]\ar[r]^(.4){\cong}&IH^*(U)\oplus
F^*(U)\ar[d]\ar[r]&IH^*(U)\ar[d] \\
V^*_{\cN_U}\ar@{^(->}[r]&V^*_{\hat{\cN}_U}\ar[r]&IH^*(\hat{\cN}_U\git
K )\ar[r]^(0.4){\cong}&IH^*(\cN_U\git K)\oplus
F^*(U)\ar[r]&IH^*(\cN_U\git G)}$$ where $\hat{\cN}_U$ is the
preimage of $\cN_U$ in $\hatm^{ss}$. Since $\mathcal{F}^{\cdot}$
is supported over $X_{(\lih)}=GM^{ss}_H\git G$, the vertical map
for $F^*(U)$ is the identity map.

Now let $\zeta$ be a nonzero element in $B^*(U)$. We know $\zeta$
is mapped to a nonzero element, say $\eta$ in $F^*(U)$. In the
above diagram, $\eta$ is mapped to zero in $IH^*(\cN_U\git G)$.
Then by Proposition \ref{lem8.6a}, $\zeta|_{\cN_U}=0$ and thus
$\eta=0$. This is a contradiction! So we proved that $B^*(U)=0$.
\end{proof}

It remains to prove Proposition \ref{lem8.6a}. This is a
consequence of the next three lemmas.

Let $x$ and $\cN_x$ be as in Definition \ref{5.3}. By \cite{k2}
Corollary 5.6 and \cite{k5} Lemma 1.21, we have an isomorphism
\begin{equation}\lab{eqisomNu}
H^*_K(\cN_U)\cong \left[H^*\left(M^{ss}_H\cap \phi^{-1}(U)\git
N^H_0\right)\otimes H^*_H(\cN_x)\right]^{\pi_0N^H}\end{equation}
 from a degenerating spectral sequence.

\begin{lemma}\lab{lemlast1}
Via the isomorphism (\ref{eqisomNu}), we have
\begin{equation}\lab{eqisomNw}
V^*_{\cN_U}\cong \left[H^*\left(M^{ss}_H\cap \phi^{-1}(U)\git
N^H_0\right)\otimes V^*_{\cN_x}\right]^{\pi_0N^H}.\end{equation}
\end{lemma}

\begin{proof}
Let $L\in \mathcal{R}(M)$.  If $L$ is not conjugate to a subgroup
of $H$, there is no $L$-fixed point in $\cN$. Hence, after
conjugation if necessary, we may assume that $L\subset P$. Let
$\cN_{U,L}$ be the $L$-fixed subset of $\cN_U$. For $V^*_{\cN_U}$
we have to consider the map
\begin{equation}\lab{eqforp8.5.1} H^*_K(\cN_U)\to
H^*_K(K\times_{N^L}\cN_{U,L})\cong H^*_{N^L}(\cN_{U,L})\to
[H^*_{N^L_0/L}(\cN_{U,L})\otimes H^{\ge n_L} _L]^{\pi_0N^L}.
\end{equation}

It is obvious that $\cN_{U,L}$ is a vector bundle over
$M^{ss}_L\cap GM^{ss}_H\cap \phi^{-1}(U)$. Using Lemma
\ref{lem8.1a}, it is easy to check that there are $k_1,\cdots,
k_s$ in $K$ so that
$$M^{ss}_L\cap GM^{ss}_H\cap
\phi^{-1}(U)=\bigsqcup_iN^L_0M^{ss}_{H_i}\cap \phi^{-1}(U)$$ where
$H_i=k_iHk_i^{-1}$. Since $H$ is maximal in $\mathcal{R}(M)$ we
have isomorphisms $$N^L_0M^{ss}_{H_i}\cap \phi^{-1}(U)\cong
N^L_0\times_{N^L_0 \cap N^{H_i}}M^{ss}_{H_i}\cap \phi^{-1}(U)\cong
N^L_0/L\times_{N^L_0 \cap N^{H_i}/L}M^{ss}_{H_i}\cap
\phi^{-1}(U)$$ by \cite{k2} 5.6 and hence we have
\begin{equation}\lab{eqp8.5.2}
H^*_{N^L_0/L}(\cN_{U,L})\cong \bigoplus_iH^*_{N^L_0\cap
N^{H_i}/L}\left(M^{ss}_{H_i}\cap \phi^{-1}(U)\right).
\end{equation}

If we apply Lemma \ref{lem8.1a} with $L\subset H$ as subgroups of
$N^H$, we deduce that there exist $g_1,\cdots, g_t$ in $N^H$ such
that $$\{k\in N^H\,|\, k^{-1}Lk\subset H\}=\bigsqcup_j (N^H\cap
N^L_0)g_jH.$$ But for $k\in N^H$, $k^{-1}Lk\subset k^{-1}Hk=H$ and
hence we have
$$N^H=\bigsqcup_j
(N^H\cap N^L_0)g_jH.$$ This implies that the natural embedding
$$N^H\cap N^L_0/H\cap N^L_0\hookrightarrow N^H/H$$ is of finite
index. In particular, the identity component of $N^H\cap
N^L_0/H\cap N^L_0$ is naturally isomorphic to the identity
component $N^H_0/H$ of $N^H/H$.

Let $N_i$ be the identity component of $N^{H_i}\cap N^L_0$ and put
$S_i=N_i\cap H_i$. Then $N_i/S_i\cong N^{H_i}_0/H_i$. Therefore,
$$H^*_{N^L_0\cap N^{H_i}/L}\left(M^{ss}_{H_i}\cap
\phi^{-1}(U)\right)$$ is the $\pi_0(N^{H_i}\cap N^L_0)$-invariant
part of
\begin{equation}\lab{eq8.com}
\begin{array}{ll}
H^*_{N_i/L}\left(M^{ss}_{H_i}\cap \phi^{-1}(U)\right)&\cong
H^*_{N_i/S_i}\left(M^{ss}_{H_i}\cap \phi^{-1}(U)\right)\otimes
H^*_{S_i/L}\\
& \cong H^*_{N^{H_i}_0/H_i} \left(M^{ss}_{H_i}\cap
\phi^{-1}(U)\right)\otimes H^*_{S_i/L}\\
& \cong H^*\left(M^{ss}_{H_i}\cap \phi^{-1}(U)\git
N^{H_i}_0\right)\otimes H^*_{S_i/L}\end{array}
\end{equation}
The first isomorphism in (\ref{eq8.com}) came from \cite{k5} Lemma
1.21. Our interest lies in finding the kernel of
(\ref{eqforp8.5.1}). Combining (\ref{eqforp8.5.1}),
(\ref{eqp8.5.2}) and (\ref{eq8.com}), we see that $V^*_{\cN_U}$ is
the intersection of the kernels of
\begin{equation}\lab{eqp8.5.3}
H^*_K(\cN_U)\to \bigoplus_i H^*\left(M^{ss}_{H_i}\cap
\phi^{-1}(U)\git N^{H_i}_0\right)\otimes H^*_{S_i/L}\otimes H^{\ge
n_L} _L \end{equation}
for all $L\in\mathcal{R}(M)$.

Now observe that the spaces that appear in (\ref{eqforp8.5.1}) lie
over $GM^{ss}_H\git G\cong M^{ss}_H\git N^H$. Applying spectral
sequence, we get a homomorphism of spectral sequences whose
$E_2$-terms give us \small
\begin{equation}\lab{eq8.19c}
\left[H^*\left(M^{ss}_H\cap \phi^{-1}(U)\git N^H_0\right)\otimes
H^*_H(\cN_x)\right]^{\pi_0N^H}\to \bigoplus_i
H^*\left(M^{ss}_{H_i}\cap \phi^{-1}(U)\git N^{H_i}_0\right)\otimes
H^*_{S_i/L}(\cN_{x,L})\otimes H^{\ge n_L} _L.
\end{equation}\normalsize
The right side of (\ref{eq8.19c}) is isomorphic to
\begin{equation}\lab{eq8.19d}
\bigoplus_i H^*\left(M^{ss}_{H}\cap \phi^{-1}(U)\git
N^{H}_0\right)\otimes H^*_{k_i^{-1}S_ik_i/L_i}(\cN_{x,L_i})\otimes
H^{\ge n_L} _{L_i}
\end{equation}
by conjugation, where $L_i=k_i^{-1}Lk_i$. Note that
$k_i^{-1}S_ik_i/L_i$ and $H\cap N^{L_i}_0/L_i$ share the same
identity component say $S'_i$. Thus
$H^*_{k_i^{-1}S_ik_i/L_i}(\cN_{x,L_i})$ is the invariant subspace
of $H^*_{S'_i}(\cN_{x,L_i})$ with respect to a finite group
action. With the isomorphism (\ref{eq8.19d}), the fiber direction
in (\ref{eq8.19c}) is exactly the truncation map
$$
H^*_H(\cN_x)\to H^*_{S'_i}(\cN_{x,L_i})\otimes H^{\ge n_L} _{L_i}
$$
for $V^*_{\cN_x}$. Therefore taking the kernel of
(\ref{eqforp8.5.1}) for all $L\in \mathcal{R}(M)$ gives us
$$\left[H^*\left(M^{ss}_H\cap \phi^{-1}(U)\git
N^H_0\right)\otimes V^*_{\cN_x}\right]^{\pi_0N^H}.$$ So we are
done.
\end{proof}

Let $Y=\pp \left(\cN_x\oplus\cc\right)\supset \cN_x$ and $\hat{Y}$
be the blow-up of $Y$ at $0\in \cN_x$. Consider the action of $H$
on $Y$ and $\hat Y$ where $H$ acts trivially on the summand $\cc$.

\begin{lemma}\lab{lemlast2}
The actions of $H$ on $Y$ and $\hat Y$ are weakly balanced.
\end{lemma}
\begin{proof}
We identify the hyperplane at infinity $Y-\cN_x$ with $\pp\cN_x$.
The equation for the hyperplane at infinity is $H$-invariant and
hence $Y^{ss}$ contains $\cN_x$. Therefore $Y^{ss}=\cN_x\cup
\pp\cN_x^{ss}$. Let $0\neq y\in\cN_x$ and suppose be the
corresponding point $\overline{y}$ in $\pp\cN_x$ is semistable. By
\cite{k2} Lemma 4.3, $y$ is fixed by $L$ if and only if
$\overline{y}$ is fixed by $L$. Hence we may consider only points
in $\cN_x$.

Let $\mathrm{Stab}(x)=P$. Since we could use any $x\in\mu^{-1}(0)$
as long as the infinitesimal stabilizer is $\mathrm{Lie}\,H$, we
assume that $P$ is minimal among the stabilizers whose Lie algebra
is $\mathrm{Lie}\,H$ so that $\cN_x=W$ in the notation of Lemma
\ref{local}.

First since the action of $K$ on $M$ is weakly balanced, the
action of $H$ on $\cN_x$ is weakly linearly balanced and so is the
action of $H\cap N^L/L$ on the $L$-fixed subspace $\cN_{x,L}$.
Hence we checked the weakly balanced condition for $H$.

Now let $0\neq y\in \mu_{\cN_x}^{-1}(0)$ and let $L$ be the
identity component of the stabilizer of $y$. Then from
(\ref{eq8.6}), we see that the normal space to $GM^{ss}_L$ in
$M^{ss}$ is the same as the normal space to $P_{\cc}W^L$ in $W$ at
a generic point. But $H$ is the identity component of $P$ and
$\cN_x=W$. Hence the normal space to $GM^{ss}_L$ in $M^{ss}$ at a
generic point is isomorphic to the normal space to
$H_{\cc}\cN_{x,L}$ in $\cN_{x}$ at a generic point. Moreover since
the diffeomorphism in Lemma \ref{local} is equivariant, the
actions of $L$ on the normal spaces are identical. According to
Definition \ref{5.3}, the weakly balanced condition is purely
about the action of $L$ on the normal spaces to $H_{\cc}\cN_{x,L}$
for all $L$. Because the action of $K$ on $M$ is weakly balanced,
we deduce that the action of $H$ on $\cN_x$ is weakly balanced. So
we proved that the action of $H$ on $Y$ is weakly balanced.

Note that $\hat Y$ is the first blow-up in the partial
desingularization for $Y$ and hence $r(\hat Y)<r=r(M)$ by
\cite{k2} 6.1. By Remark \ref{remWBhat} (2), the action on $\hat
Y$ is also weakly balanced.
\end{proof}

If we let $\cN_x=\mathrm{Spec}\,A$ for some polynomial ring $A$,
then by definition \cite{mfk} there are homeomorphisms
$$\cN_x\git H\cong \mathrm{Spec}\,A^H\qquad
\pp\cN_x\git H\cong \mathrm{Proj}\, A^H$$ i.e. $\cN_x\git H$ is
the affine cone of $\pp\cN_x\git H$ and $Y\git H$ is the
projective cone of $\pp\cN_x\git H$.

Since the equation for the hyperplane at infinity in $Y$ is
invariant, by construction of $Y\git G$ (\cite{mfk} Chapter 1
\S4), it is obvious that the preimage of $\cN_x\git H$ by the GIT
quotient map $\phi_Y:Y^{ss}\to Y\git H$ is $\cN_x$. Because the
action of $H$ on $Y$ is weakly balanced by Lemma \ref{lemlast2},
the Kirwan map for $Y$ restricted to $\cN_x$
$$\kappa^{ss}_{\cN_x}:H^*_H(\cN_x)\to
IH^*(\cN_x\git H)$$ is surjective by Proposition \ref{prop6.2b}.

Now the (sheaf-theoretic) Kirwan map gives us a homomorphism of
the spectral sequence for $H^*_K(\cN_U)$ to the spectral sequence
for $IH^*(\cN_U\git K)$. At $E_2$-level, we have a homomorphism
\small
$$
\left[H^*\left(M^{ss}_H\cap \phi^{-1}(U)\git N^H_0\right)\otimes
H^*_H(\cN_x)\right]^{\pi_0N^H}\to \left[H^*\left(M^{ss}_H\cap
\phi^{-1}(U)\git N^H_0\right)\otimes IH^*(\cN_x\git H
)\right]^{\pi_0N^H}.$$ \normalsize
 But we know the left side
spectral sequence degenerates (\ref{eqisomNu}) and $H^*_H\cong
H^*_H(\cN_x)$ surjects onto $IH^*(\cN_x\git H)$. Therefore the
right side spectral sequence also degenerates and we have an
isomorphism (\cite{k3} p495)
\begin{equation} \lab{eqisomNv}
IH^*(\cN_U\git K)\cong \left[H^*\left(M^{ss}_H\cap
\phi^{-1}(U)\git N^H_0\right)\otimes IH^*(\cN_x\git H
)\right]^{\pi_0N^H}
\end{equation}

In view of (\ref{eqisomNw}) and (\ref{eqisomNv}), the proof of
Proposition \ref{lem8.6a} is complete if we show the following
lemma.
\begin{lemma}
$V^*_{\cN_x}\cong IH^*(\cN_x\git H)$.
\end{lemma}
\begin{proof}
As in the proof of Lemma \ref{lemlast2}, we may assume that
$\mathrm{Stab}(x)$ is minimal among the stabilizers of points in
$\mu^{-1}(0)$ whose Lie algebra is $\mathrm{Lie}\, H$. So we may
use Lemma \ref{local} with $\cN_x=W$.

As $\cN_x\git H$ is a cone with vertex point $*$, it is well-known
that
$$IH^i(\cN_x\git H)=0$$ for $i\ge n_H$ and if $i<n_H$, we have
$$IH^i(\cN_x\git H)\cong IH^i(\cN_x\git H-*) \cong IH^i(\Cal
N_x-\phi^{-1}_x(*)\git H)$$ where $\phi_x:\cN_x\rightarrow
\cN_x\git H$ is the GIT quotient map and the following diagram
commutes:
\begin{equation}\lab{D4.6}
\begin{CD}
IH^i(\cN_x\git H)  @>>>   IH^i(\cN_x-\phi^{-1}_x(*)\git H)\\
@V{\phi_{x,H}^*}VV   @V{\phi_{x,H}^*}VV\\
H^i_H(\cN_x)  @>>>    H^i_H(\cN_x-\phi^{-1}_x(*))\\
\end{CD}
\end{equation}
As $\phi^{-1}_x(*)$ is the union of the complex cones over the
unstable strata of $\pp\cN_x$ and the real codimension of each
unstable stratum is greater than $n_H$ by the weakly balanced
condition, the real codimension of $\phi^{-1}_x(*)$ is greater
than $n_H$. Hence, the bottom horizontal map is an isomorphism.

We claim that the the right vertical in (\ref{D4.6}) is injective
and the image is $V^i_{\cN_x-\phi_x^{-1}(*)}$. Consider the
following commutative diagram
$$\xymatrix{
\hat{Y}^{ss}\ar[r]^{\delta}\ar[d]_{\tilde{\gamma}}& Y^{ss}\ar[d]^{\gamma}\\
\hat{Y}\git H\ar[r]^{\epsilon}&Y\git H }$$ we see that
$$\tilde{\gamma}^{-1}\epsilon^{-1}(\cN_x\git
H-*)=\delta^{-1}\gamma^{-1}(\cN_x\git
H-*)=\delta^{-1}\left(\cN_x-\phi_x^{-1}(*)\right)\cong
\cN_x-\phi_x^{-1}(*).$$ For the last isomorphism, observe that
$\cN_x-\phi_x^{-1}(*)$ does not intersect with the blow-up center.
Moreover, if $y\in \cN_x-\phi_x^{-1}(*)$, the closure of
$H_{\cc}\cdot y$ does not meet $0$ (if it does, the point belongs
to $\phi_x^{-1}(*)$) and thus $y$ is semistable as a point in
$\hat{Y}$ by \cite{k2} Remark 7.7.

Because the action of $H$ on $\hat Y$ is weakly balanced and
$r(\hat Y)<r=r(M)$, Theorem \ref{splittingtheorem} is true for
$\hat Y$. In particular, if we apply the theorem for the preimage
$\cN_x-\phi_x^{-1}(*)$ of the open set $\epsilon^{-1}(\cN_x\git
H-*)$ by $\tilde{\gamma}$ we obtain the isomorphism
$$IH^i(\cN_x-\phi^{-1}_x(*)\git H)\cong V^i_{\cN_x-\phi_x^{-1}(*)}.$$

Therefore, from (\ref{D4.6}), it suffices to show that
$$V^i_{\cN_x}\cong V_{\cN_x-\phi^{-1}_x(*)}^i$$ for $i<n_H$ by
restriction. To see this, we note once again that for $i<n_H$,
$$H^i_H(\cN_x)\rightarrow H^i_H(\cN_x-\phi^{-1}_x(*))$$ is an
isomorphism and the same is true for
$$H^j_{H\cap N_0^L/L}(\mathcal
N_{x,L}) \otimes H^k_{L}\rightarrow H^j_{H\cap N_0^L/L}(\mathcal
N_{x,L}-\phi^{-1}_x(*)) \otimes H^k_{L}$$ if $k\ge n_L$,
$j<n_H-k\le n_H-n_L$,  because the action of $H\cap N_0^L/L$ on
$\cN_{x,L}$ is weakly linearly balanced.
\end{proof}


\section{Examples}\lab{section9}

 Let $$P^K_t(W)=\sum_{i\ge 0} t^i \dim  H^i_K(W)$$
$$IP_t(W)=\sum_{i\ge 0} t^i \dim  IH^i(W)$$
be the Poincar\'e series.

\subsection{$\Bbb C^*$-action on projective space}
Consider a $\Bbb C^*$-action on $M=\Bbb P^n$ via a representation
$\Bbb C^*\rightarrow GL(n+1)$. Let $n_+$, $n_0$, $n_-$ be the number
of positive, zero, negative weights.  Suppose
the action is weakly balanced, i.e. $n_+=n_-$.
In this case, we can easily compute the intersection Betti nubmers by
Theorem \ref{splittingtheorem}.

From the equivariant Morse theory \cite{k1},
\begin{equation}
\begin{split}
P_t^{S^1}((\Bbb P^n)^{ss})&=
\frac{1+t^2+\cdots +t^{2n}}{1-t^2}-\frac{t^{2n_0+2n_+}+\cdots +
t^{2n}}{1-t^2}-\frac{t^{2n_0+2n_-}+\cdots +
t^{2n}}{1-t^2}\\
&=\frac{1+t^2+\cdots +t^{2n_0+2n_--2}-t^{2n_0+2n_+}-\cdots
-t^{2n}}{1-t^2}.\end{split}\end{equation}
In this case, $\Cal R=\{S^1\}$ and $M^{ss}_{S^1}=\Bbb P^{n_0-1}$,
$n_H=n_++n_--1=2n_--1$.

As $H^*_{S^1}(\Bbb P^n)\rightarrow H^*_{S^1}((\Bbb P^n)^{ss})
\rightarrow H^*_{S^1}(\Bbb P^{n_0-1})$ is surjective, we
have only to subtract out the Poincar\'e series of
$\oplus_{i\ge 2n_-}H^*(\Bbb P^{n_0-1})\otimes H^i_{S^1}$,
which is precisely
$$\frac{t^{2n_-}(1+t^2+\cdots +t^{2n_0-2})}{1-t^2}.$$
Therefore,
$$IP_t(\Bbb P^n\git \Bbb C^*)=\frac{1+t^2+\cdots +t^{2n_--2}-
t^{2n_0+2n_+}-t^{2n_0+2n_++2}\cdots -t^{2n}}{1-t^2}$$
which is a palindromic polynomial of degree $2n-2$.

\subsection{Ordered $2n$-tuples of points of $\Bbb P^1$}
Let us consider $G=SL(2)$ action on the set $M=(\Bbb P^1)^{2n}$
of ordered $2n$-tuples of points in $\Bbb P^1$ as M\"obius transformations.
Then the semistable $2n$-tuples are those containing no point of $\Bbb P^1$
strictly more than $n$ times and the stable points are those containing no
point at least $n$ times. Let $H$ be the maximal torus of $G$. Then
$$\Cal R(M)=\{H\}$$
$$M^{ss}_H=\{q_I | I\subset \{1,2,3,\cdots, 2n\}, |I|=n\}$$
where $q_I\in M$, the j-th component of which is $\infty$ if $j\in I$, $0$
otherwise. Therefore, the action is weakly balanced.

The normalizer $N^H$ satisfies $N^H/H=\Bbb Z/2$ and $N^H_0=H$. Hence,
$$H^*_{N^H}(M^{ss}_H)=[\oplus_{|I|=n}H^*_{H}]^{\Bbb Z/2}.$$
The $\Bbb Z/2$ action interchanges $\infty$ and $0$, i.e. $q_I$ and $q_{I^c}$,
and therefore, from now on,
we only think of those $I$'s that contain $1$ so that we can
forget the $\Bbb Z/2$ action.

$H^*_K(M)=H^*_K((\Bbb P^1)^{2n})$
has generators $\xi_1,\xi_2,\cdots ,\xi_{2n}$ of
degree 2 and $\rho^2$ of degree 4, subject to the relations $\xi_j^2=\rho^2$
for $1\le j\le 2n$. The $I$-th component of the restriction map
$$H^*_K(M)\rightarrow H^*_{N^H}(M^{ss}_H)=\oplus_IH^*_{H}$$
maps $\rho^2$ to $\rho^2$ and $\xi_j$ to $\rho$ if $j\in I$, $-\rho$ otherwise.
From \cite{k1},
$$P^K_t(M^{ss})=
\frac{(1+t^2)^{2n}}{1-t^4}-\sum_{n<r\le 2n}\binom{2n}r\frac{t^{2(r-1)}}{1-t^2}.
$$

\begin{proposition} $$
IP_t(M\git G)=P^K_t(M^{ss})-\frac12\binom{2n}n\frac{t^{2n-2}}{1-t^2}
$$\end{proposition}
\begin{proof} By Theorem \ref{splittingtheorem}, we have only to subtract out
the Poincar\'e series of
$$Im\{H^*_K(M^{ss})\rightarrow \oplus_I H^*_{H}\}\,\,\,
\cap \,\,\, \{\oplus_I\oplus_{i\ge n_H} H^i_{H}\}$$
where $n_H$ is in this case $2n-3$. By the lemma below, which is essentially
combinatorial, the image contains $\oplus_I\oplus_{i\ge 2n-3}
H^i_{H}$ and thus the intersection is $\oplus_I\oplus_{i>2n-3}
H^i_{H}$, whose Poincar\'e series is precisely
$$\frac12\binom{2n}n\frac{t^{2n-2}}{1-t^2}.$$
So we are done. \,\,\,\,\,\,$\square$\end{proof}

\begin{lemma} The restriction map
$H^{2k}_K(M^{ss})\rightarrow \oplus_I H^{2k}_{H}$
is surjective for $k\ge n-1$.
\end{lemma}
\begin{proof} It is equivalent to show that
$H^{2k}_K(M)\rightarrow \oplus_I H^{2k}_{H}$ is surjective. Let
$\xi=\xi_2+\cdots +\xi_{2n}$ and consider, for each $I=(1,i_2,\cdots,i_n)$,
$$\eta_I=(\xi-\xi_{i_2})^{k_2}\cdots (\xi-\xi_{i_n})^{k_n}.$$
Then since $\xi|_{q_J}=-\rho$ for all $J$, $\eta_I|_{q_J}=(-2\rho)^k$ if $J=I$
and $0$ otherwise, where $k=k_2+\cdots +k_n$, $k_i\ge 1$.
Therefore, the images of those $\eta_I$ span $\oplus_I H^{2k}_{H}$ for any
$k\ge n-1$ and thus the restriction is surjective.\,\,\,\,\,\,$\square$\end{proof}

\subsection{Intersection pairing}

Consider the $\Bbb C^*$ action on $\Bbb P^7$ by
a representation with weights $+1,0,-1$ of multiplicity $3,2,3$ respectively.
Hence, $n_+=n_-=3$, $n_0=2$. Then
$$H^*_{S^1}(\Bbb P^7)=\Bbb C[\xi,\rho]/\langle
\xi^2(\xi-\rho)^3(\xi+\rho)^3\rangle$$
where $\xi$ is a generator in $H^2(\Bbb P^7)$ and $\rho$ is a generator
in $H^2_{S^1}$.

The equivariant Euler classes for the two unstable strata are
$\xi^2(\xi-\rho)^3$, $\xi^2(\xi+\rho)^3$ respectively. Therefore,
$$H^*_{S^1}((\Bbb P^7)^{ss})
=\Bbb C[\xi,\rho]/\langle\xi^2(\xi-\rho)^3,\, \xi^2(\xi+\rho)^3\rangle $$
A Gr\"obner basis for the relation ideal is
$$\{\xi^5+3\xi^3\rho^2, \xi^4\rho+\frac13\xi^2\rho^3, \xi^3\rho^3,
\xi^2\rho^5\}$$
where $\xi>\rho$. Hence as a vector space,
$$H^*_{S^1}((\Bbb P^7)^{ss})
=\Bbb C\{\xi^i\rho^j | i=0,1,\, j\ge 0\}\oplus
\Bbb C\{\xi^i\rho^j | 2i+j<9, i\ge 2, j\ge 0\}$$
By definition, as $n_{S^1}=5$,
we remove $\Bbb C\{\xi^i\rho^j | i=0,1,\, j\ge 3\}$ to get
$$V=\oplus_{0\le i\le 6} V^{2i}$$
$$V^0=\Bbb C,\,\,\,\,\, V^2=\Bbb C\{\rho,\xi\},\,\,\,\,\,\,
V^4=\Bbb C\{\rho^2,\xi\rho, \xi^2\},$$
$$V^6=\Bbb C\{\xi\rho^2,\xi^2\rho,\xi^3\},\,\,\,\,\,
V^8=\Bbb C\{\xi^2\rho^2,\xi^3\rho,\xi^4\},$$
$$V^{10}=\Bbb C\{\xi^2\rho^3,\xi^3\rho^2\},\,\,\,\,\,\,
V^{12}=\Bbb C\{\xi^2\rho^4\}.$$
First, consider the pairing $V^2\otimes V^{10}\rightarrow V^{12}$.
As $\rho (\xi^2\rho^3)=\xi^2\rho^4$, $\rho(\xi^3\rho^2)=\xi^3\rho^3=0$,
$\xi(\xi^2\rho^3)=\xi^3\rho^3=0$, $\xi(\xi^3\rho^2)=\xi^4\rho^2=-\frac13
\xi^2\rho^4$, the pairing matrix is up to a constant
$$\left(\begin{matrix} 1 &0\\0&-\frac13\end{matrix}\right)$$
The determinant is $-\frac13\ne 0$ and the signature is $0$.

Next, consider the pairing $V^4\otimes V^8\rightarrow V^{12}$.
One can similarly use the Gr\"obner basis to compute the pairing
as above. The pairing matrix is up to a constant
$$\left(\begin{matrix} 1&0&-\frac13\\ 0&-\frac13&0\\-\frac13&0&1\end{matrix}
\right)$$
The determinant is $-\frac{8}{27}\ne0$ and the signature is 1.

Similarly, the intersection pairing matrix for
$V^6\otimes V^6\rightarrow V^{12}$ is up to a constant
$$\left(\begin{matrix} 1&0&-\frac13\\ 0&-\frac13&0\\-\frac13&0&1\end{matrix}
\right)$$
The determinant is $-\frac{8}{27}\ne0$ and the signature is 1.

In this way, one can compute the intersection pairing for any $n_0, n_-=n_+$.

\vspace{1cm}
In a subsequent paper, we will compute the intersection pairing of the
moduli spaces of holomorphic vector bundles over a Riemann surface
 of any rank and degree, using the nonabelian localization theorem
of Jeffrey and Kirwan.

\newcommand{\etalchar}[1]{$^{#1}$}


\begin{thebibliography}{MFK94}

\bibitem[AB82]{ab2}
M.F. Atiyah and R.~Bott.
\newblock The {Y}ang--{M}ills equations over {R}iemann surfaces.
\newblock {\em Phil. Trans. Roy. Soc. Lond.}, A308:532--615, 1982.

\bibitem[B{\etalchar{+}}84]{bo}
A.~Borel et~al.
\newblock {\em Intersection cohomology}.
\newblock Number~50 in Progress in mathematics. Birkh\"auser, 1984.

\bibitem[BBD82]{BBD}
A.~Beilinson, J.~Bernstein, and P.~Deligne.
\newblock Faisceaux pervers.
\newblock {\em Ast\'erisque}, 100, 1982.
\newblock Proc. C.I.R.M. conf\'erence: Analyse et topologie sur les espaces
  singuliers.

\bibitem[BL94]{BL}
J.~Bernstein and V.~Lunts.
\newblock {\em Equivariant sheaves and functors}.
\newblock Lecture Notes in Math. 1578, Springer, 1994.


\bibitem[GeMa]{GeM}
S.~Gelfand and Y.~Manin.
\newblock {\em Methods of homological algebra}.
\newblock Springer-Verlag, 1996.

\bibitem[GH94]{GH}
P.~Griffiths and J.~Harris.
\newblock {\em Principles of Algebraic Geometry}.
\newblock Wiley, 1994.

\bibitem[GM80]{GM1}
M.~Goresky and R.~MacPherson.
\newblock Intersection homology theory {I}.
\newblock {\em Topology}, 19:135--162, 1980.

\bibitem[GM83]{GM2}
M.~Goresky and R.~MacPherson.
\newblock Intersection homology theory {II}.
\newblock {\em Inventiones Mathematicae}, 71:77--129, 1983.

\bibitem[GM85]{gm6}
M.~Goresky and R.~MacPherson.
\newblock Lefschetz fixed point theorem for intersection homology.
\newblock {\em Comment. Math. Helvetici}, 60:366--391, 1985.

\bibitem[Jef94]{j}
L.C. Jeffrey.
\newblock Extended moduli spaces of flat connections on {R}iemann surfaces.
\newblock {\em Math. Annalen}, 298:667--692, 1994.

\bibitem[JK98]{JK98}
L.C. Jeffrey and F.C. Kirwan.
\newblock Intersection theory on moduli spaces of holomorphic
bundles of arbitrary rank on a Riemann surface.
\newblock {\em Ann. Math.}, 148:109--196, 1998.

\bibitem[Kie]{kiem2}
Y.-H. Kiem.
\newblock Intersection cohomology of representation spaces of surface groups.
\newblock Preprint.

\bibitem[Kir84]{k1}
F.~Kirwan.
\newblock {\em Cohomology of Quotients in Symplectic and Algebraic Geometry}.
\newblock Number~34 in Mathematical Notes. Princeton University Press, 1984.

\bibitem[Kir85]{k2}
F.~Kirwan.
\newblock Partial desingularisations of quotients of nonsingular varieties and
  their {B}etti numbers.
\newblock {\em Annals of Mathematics}, 122:41--85, 1985.

\bibitem[Kir86a]{k5}
F.~Kirwan.
\newblock On the homology of compactifications of moduli spaces of vector
  bundles over a riemann surface.
\newblock {\em Proc. Lon. Math. Soc.}, 53:237--266, 1986.

\bibitem[Kir86b]{k3}
F.~Kirwan.
\newblock Rational intersection homology of quotient varieties.
\newblock {\em Inventiones Mathematicae}, 86:471--505, 1986.

\bibitem[Kir92]{k6}
F.~Kirwan.
\newblock The cohomology rings of moduli spaces of bundles over riemann
  surfaces.
\newblock {\em Jour. A.M.S.}, 5:853--906, 1992.

\bibitem[Kir94]{k7}
F.~Kirwan.
\newblock Geometric invariant theory and the {A}tiyah-{J}ones conjecture.
\newblock {\em Proc. S. Lie Mem. Conf.}, pages 161--188, 1994.

\bibitem[KW]{kw}
Y.-H. Kiem and J.~Woolf.
\newblock The cosupport axiom, equivariant cohomology and the intersection
  cohomology of certain symplectic quotients.
\newblock Preprint.

\bibitem[MFK94]{mfk}
D.~Mumford, J.~Fogarty, and F.~Kirwan.
\newblock {\em Geometric invariant theory}.
\newblock Springer--Verlag, third edition, 1994.

\bibitem[MS99]{mes}
E.~Meinrenken and R.~Sjamaar.
\newblock Singular reduction and quantization.
\newblock {\em Topology}, 38(4):699--762, 1999.

\bibitem[New78]{New}
P.~Newstead.
\newblock {\em Introduction to moduli problems and orbit spaces}.
\newblock Tata institute of fundamental research, Bombay, 1978.

\bibitem[Sja95]{Sja}
R.~Sjamaar.
\newblock Holomorphic slices, symplectic reduction and multiplicities of
  representations.
\newblock {\em Ann. Math.}, 141:87--129, 1995.

\bibitem[SL91]{SL}
R.~Sjamaar and E.~Lerman.
\newblock Stratified symplectic spaces and reduction.
\newblock {\em Annals of Maths}, 134:375--422, 1991.

\end{thebibliography}
\end{document}